\documentclass[11 pt, psamsfonts, leqno]{amsart}
\usepackage{amssymb}
\usepackage{amsmath}	

 \calclayout
\def\ignore#1{\relax}

\usepackage{bm}

\ignore{
\usepackage[scaled=0.92]{helvet}	
\usepackage[subscriptcorrection,slantedGreek,nofontinfo,  mtphbi]{mtpro2}
}

\usepackage[plainpages=false, pdfpagelabels,  colorlinks=true, pdfstartview=FitV, linkcolor=blue, citecolor=blue, urlcolor=blue]{hyperref}

\usepackage{epic}
\usepackage{eepicemu}

\usepackage{graphicx}
\usepackage{epstopdf}

\usepackage{paralist}

\usepackage{maxidiagram}

\newcommand\inlinegraphic[2][{scale=1.0}]{\begin{array}{c} \includegraphics[#1]{./EPS/#2}\end{array}}

\numberwithin{equation}{section}
\numberwithin{figure}{section}

\def\Z{{\mathbb Z}}

\def\Q{{\mathbb Q}}

\def\eps{\varepsilon}
 
\def\Ad{{\rm Ad}}
\def\End{{\rm End}}
\def\Hom{{\rm Hom}}

\def\powerpm{^{\pm 1}} 
\def\u #1 #2{\mathcal U(#1, #2)}  
\def\uhat #1 #2{\widehat{\mathcal U}(#1, #2)}  

\def\varLambdaHat{{\widehat{\varLambda}}}

\def\kt #1{{KT_{#1}}}  
\def\akt #1{\widehat{KT}_{#1}}  
\def\abmw #1{\widehat{W}_{#1}}  
\def\bmw #1{W_{#1}}  
\def\braid #1{\mathcal B_{#1}} 
\def\affbraid  #1{\widehat{\mathcal B}_{#1}}
\def\ahec #1{\widehat{H}_{#1}}  
\def\hec #1{H_{#1}}
\def\w #1 #2{\bmw {#1}^{(#2)}}  
\def\V #1 #2{V_{#1}^{(#2)}}  
\def\k #1 #2{ KT_{#1}^{(#2)}}
 
\def\p #1{\bm {#1}}
\def\pbar #1{\bar{\p #1}}

\def\bdry{\partial}

\def\inv{^{-1}}
\def\la{\lambda}

\def\ubold{{\bm u}}

\def\id{{\rm id}}

\def\rhobold{{\bm \rho}} 
\def\deltabold{{\bm \delta}}

\def\qbold{{\bm q}}

\def\boldu{{\bm u}}

\def\A{\mathbb A}
\def\B{\mathbb B}

\def\hods{\unskip\kern.55em\ignorespaces}

\theoremstyle{plain}
\newtheorem{theorem}{Theorem}[section]

\theoremstyle{plain}
\newtheorem{proposition}[theorem]{Proposition}

\theoremstyle{plain}
\newtheorem{corollary}[theorem]{Corollary}

\theoremstyle{plain}
\newtheorem{lemma}[theorem]{Lemma}

\theoremstyle{definition}
\newtheorem{definition}[theorem]{Definition}
\theoremstyle{definition}
\newtheorem{example}[theorem]{Example}

\theoremstyle{definition}
\newtheorem{remark}[theorem]{Remark}

\theoremstyle{remark}

\title[Cyclotomic BMW algebras]{Cyclotomic Birman--Wenzl--Murakami algebras,  I:\\ Freeness and realization as tangle algebras }

\author{Frederick M. Goodman}
\address{ Department of Mathematics\\ University of Iowa\\ Iowa
City, Iowa}
\email{ goodman@math.uiowa.edu} 
\author{Holly Hauschild Mosley}
\address{Department of Mathematics\\ Grinnell College \\
Grinnell, Iowa}
\email{HAUSCHIL@GRINNELL.EDU}


\begin{document}
 \baselineskip=16pt 
 \maketitle

March, 2008

\begin{abstract}  The cyclotomic Birman-Wenzl-Murakami algebras are quotients
of the affine BMW algebras in which the affine generator satisfies a polynomial relation.  We show that the cyclotomic BMW algebras are free modules over any admissible, integral ground ring,  and that they are isomorphic to cyclotomic versions of the Kauffman tangle algebras.
\end{abstract}

\setcounter{tocdepth}{1}
\tableofcontents
\section{Introduction}  

This paper and the companion paper ~\cite{GH3} continue the study of affine and cyclotomic Birman--Wenzl-- Murakami  (BMW) algebras, which we began in  ~\cite{GH1}.  

\subsection{Background}  \label{subsection: introduction background}
The origin of the BMW algebras was in knot theory.  Kauffman defined ~\cite{Kauffman} an invariant of regular isotopy for  links in $S^3$, determined by certain skein relations.
 Birman and Wenzl ~\cite{Birman-Wenzl} and independently Murakami ~\cite{Murakami-BMW}  then defined  a family of quotients of the braid group algebras, and showed that Kauffman's invariant could be recovered from a trace on these algebras.  These (BMW) algebras were  defined by generators and relations, but were implicitly modeled on certain algebras of tangles whose definition  was subsequently made explicit by Morton and Traczyk ~\cite{Morton-Traczyk}, as follows:
  Let $S$ be a commutative unital ring with invertible elements
$\rho$, $q$, and $\delta_0$ satisfying $\rho\inv - \rho = (q\inv -q) (\delta_0 - 1)$.  The {\em Kauffman tangle algebra}  $\kt{n, S}$  is the $S$--algebra of framed $(n, n)$--tangles in the disc cross the interval,  modulo Kauffman skein relations:
\begin{enumerate}
\item Crossing relation:
$
\quad \inlinegraphic[scale=.6]{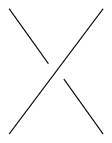} - \inlinegraphic[scale=.3]{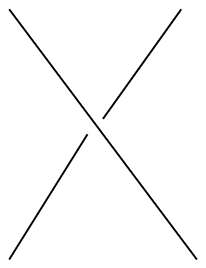} 
\quad = 
\quad
(q\inv - q)\,\left( \inlinegraphic[scale=1]{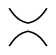} - 
\inlinegraphic[scale=1]{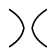}\right).
$
\item Untwisting relation:
$\quad 
\inlinegraphic{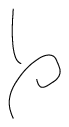} \quad = \quad \rho \quad
\inlinegraphic{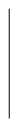} \quad\ \text{and} \quad\ 
\inlinegraphic{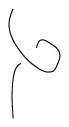} \quad = \quad \rho\inv \quad
\inlinegraphic{vertical_line}. 
$
\item  Free loop relation:  $T\, \cup \, \bigcirc = \delta_0 \, T. $
\end{enumerate}
Morton and Traczyk ~\cite{Morton-Traczyk}   showed that the $n$--strand algebra $\kt{n, S}$ is free of rank $(2n-1)!!$ as a module over $S$, and
Morton and Wassermann ~\cite{Morton-Wassermann} proved that the BMW algebras and the Kauffman tangle algebras are isomorphic.
 
It is natural to ``affinize" the BMW algebras to obtain BMW analogues of the affine Hecke   algebras of type $A$, see ~\cite{ariki-book}.  The affine Hecke algebra  can be realized geometrically as the algebra of braids in the annulus cross the interval,  modulo Hecke skein relations; this suggests defining the affine Kauffman tangle algebra $\akt{n, S}$ as the algebra of framed $(n, n)$--tangles in the annulus cross the interval, modulo Kauffman skein relations.  Turaev  ~\cite{Turaev-Kauffman-skein} showed that the resulting algebra of  $(0,0)$--tangles is a  (commutative) polynomial algebra in infinitely many variables, so it  makes sense  to absorb this polynomial algebra into the ground ring.  (The ground ring gains infinitely many parameters $\delta_j$ ($j \ge 1$)  corresponding to the generators of the polynomial algebra.)
With this,  one arrives at the definition of the affine Kauffman tangle algebra given in ~\cite{GH1},  see Definition \ref{definition affine Kauffman tangle algebra} below.
On the other hand,  H\"aring--Oldenburg ~\cite{H-O2} defined an affine version of the BMW algebras by generators and relations.  In ~\cite{GH1}, we showed that 
H\"aring--Oldenburg's affine BMW algebras are isomorphic to the affine Kauffman tangle algebras,  and we showed that these algebras are free modules over their ground ring, with a basis reminiscent of a well--known basis of affine Hecke algebras.

The affine BMW algebras arise naturally in several different contexts:

\begin{list}{$\bullet$}{
\setlength \leftmargin{12 pt}}
\item {\em  Knot theory\ Êin the solid torus.}   The Kauffman skein relations determine 
a family of invariants of links in the sold torus ~\cite{Turaev-Kauffman-skein};  the family has infinitely many parameters.  Lambropoulou ~\cite{Lambr-thesis}, ~\cite{Lambr-solid-torus}  has shown that the Jones paradigm \cite{jones-invariant} for link invariants in $S^3$ can be extended to links in the solid torus; namely, invariants of links in the solid torus can be derived from Markov traces on the braid group of the annulus (which coincides with the Artin braid group of type $B$).   The Kauffman--type invariants for links in the solid torus can be recovered from the Markov trace on affine and cyclotomic BMW algebras.  The existence and uniqueness of the Markov trace on these algebras is discussed in ~\cite{GH1} and in this paper.  The infinitely many parameters for Kauffman--type invariants enter into the definition of the affine BMW algebras;  for each choice of parameters, the unique Markov trace on the corresponding affine BMW algebra yields the Kauffman--type invariant with those parameters.
\item  {\em Quantum groups and $R$--matrices.}    The following is a brief summary of results from ~\cite{RO}.   Let $\mathfrak g$ be a complex simple Lie algebra, and $U_q \mathfrak g$  the quantum universal enveloping algebra of $\mathfrak g$.   Let $M, V$ be  finite dimensional $U_q \mathfrak g$ modules.  Let $\check R_i = \check R_{V, V}$  acting on the $i$--th and $i+1$--st tensor places in $V^{\otimes f}$, and let $\check R_0^2$ denote $\check R_{V, M} \check R_{M, V}$ action on $M \otimes V$.    Then $\sigma_i \mapsto \check R_i$  gives a representation of the $f$--string braid group in $\End_{U_q \mathfrak g}(V^{\otimes f})$.  The assignments  $\sigma_i \mapsto \check R_i$ and 
$\sigma_0^2 \mapsto \check R_0^2$ determines a representation of the braid group of type $B$ in
 $\End_{U_q \mathfrak g}(M\otimes V^{\otimes f})$.    Suppose $\mathfrak g$  is an orthogonal or symplectic Lie algebra, 
 $M$ is an irreducible representation and $V$ is the vector representation.  Then the representation of the type B braid group by $\check R$ matrices  factors through a cyclotomic BMW algebra.

\item  {\em Representations of ordinary BMW algebras.}   In the ordinary BMW algebra $W_n$,  let
$L_1 = 1$,  and $L_i  =   g_{i-1} g_{i-2} \cdots g_2 \,g_1^2\, g_2 \cdots g_{i-2} g_{i-1}$  for $2 \le i \le n$.
The elements $L_i$  are  Jucys-Murphy elements for the BMW algebras.   The representation theory of the BMW algebras is largely determined by the spectrum of the (mutually commuting) elements $L_i$  and their commutation relations with the standard generators $e_i$,  $g_i$  of the algebras.  See, for example,  ~\cite{mathas-JM},  Example 2.18,  and ~\cite{enyang-up-down}. 

For generic values of the parameters,  the BMW algebra $W_n$  is semisimple and has irreducible representations labeled by Young diagrams of size $n - 2f$  ($0 \le f \le n/2$).   The irreducible representation labeled by $\la$  has a basis indexed by up--down tableaux of length $n$ and shape $\la$, that is sequence of Young diagrams beginning with the empty diagram and ending with $\la$,  in which any two successive Young diagrams differ by the addition or deletion of a box.   Fix integers $n < N$ and Young diagrams $\la$ and $\mu$  of sizes $n-2f$  and $N - 2k$,  and consider up--down tableaux beginning with $\la$ and ending with $\mu$.   Then the algebra generated by $L_{n+1}$  and 
$e_j,  g_j$  for $n+1 \le j \le N-1$  acts on the vector space spanned by such up--down tableaux, and affords a representation of the $(N-n)$--strand affine BMW algebra.

Indeed,  one can regard a representation of the BMW algebra on a space spanned by up--down tableaux as pieced together from representations of two--strand affine BMW algebras generated by triples
$\{L_j, e_j,  g_j\}$  acting on up--down tableaux of length 2.  This point of view (applied to representations of the symmetric group)  was stressed by  Okounkov and Vershik ~\cite{okounkov-vershik-selecta},  who used it to reconstruct the representation theory of the symmetric groups {\em ab initio}.

\end{list}

\subsection{Cyclotomic BMW algebras} \label{subsection: introduction cyclotomic BMW algebras}
 In this paper and the companion paper ~\cite{GH3} we consider cyclotomic BMW algebras,  which are the BMW analogues of cyclotomic Hecke algebras ~\cite{ariki-book}.  The affine BMW algebras have a distinguished generator $x_1$, which, in the geometric (Kauffman tangle) picture is represented by a braid with one strand wrapping around the hole in the annulus cross interval.  The cyclotomic BMW algebra $\bmw{n, S, r}$ is defined to be the quotient of the affine BMW algebra $\abmw{n, S}$ in which the generator $x_1$ satisfies a monic polynomial equation
\begin{equation} \label{cyclotomic polynomial relations 0}
x_1^r + \sum_{k=0}^{r-1} a_k x_1^k = 0.
\end{equation}
with coefficients in $S$.\footnote{Actually, we will assume that  the polynomial splits in $S$.}  The cyclotomic BMW algebras were also introduced by H\"aring-Oldenburg in ~\cite{H-O2}.

In the geometric (Kauffman tangle) picture, it is more natural to convert this relation into a local skein relation:
\begin{equation} \label{cyclotomic skein relations 0}
 T_r + \sum_{k=0}^{r-1} a_k  T_k = 0,
\end{equation}
whenever $T_0, T_1, \dots, T_r$  are affine tangle diagrams that are identical in the exterior of some disc $E$ and    $T_k \cap E$  consists of one strand wrapping $k$ times around the hole in the annulus cross interval; i.e.  $T_k \cap E$ ``equals" $x_1^k$.    The cyclotomic Kauffman tangle algebra $\kt{n, S, r}$  is defined to be the quotient of the affine Kauffman tangle algebra
$\akt{n, S}$ by the cyclotomic skein relation.

A priori, the ideal  in $\akt{n, S} \cong \abmw{n,S}$ generated by the cyclotomic skein relation \ref{cyclotomic skein relations 0} is larger than the ideal generated by the polynomial relation \ref{cyclotomic polynomial relations 0}, so we have a surjective, but not evidently injective homomorphism $\varphi : \bmw{n, S, r} \to \kt{n, S, r}.$ 

Another point of view that underlines the a priori distinction between $\bmw{n, S, r}$ and
$ \kt{n, S, r}$ is the following:  Consider the affine Kauffman tangle category,  with objects   the natural numbers $0, 1, 2, \dots$ and with  $\Hom(k, \ell)$ defined as the $S$--module of 
$(k, \ell)$-- tangles in the annulus cross the interval, modulo Kauffman skein relations; thus $\akt{n, S} = \End(n)$
 in this category. Then  $\akt{*, S} = \bigoplus_{k, \ell} \Hom(k, \ell)$ is a (non--unital) algebra, containing each $\akt{n, S}$ as a subalgebra.  (Regard the elements of $\akt{*, S}$ as infinite matrices with $(k, \ell)$ entry in $\Hom(\ell , k)$.)   In $\akt{*, S}$,  let $I_*$ be the ideal generated by the polynomial relations \ref{cyclotomic polynomial relations 0}, {\em one relation for each $n$}.  (There is a  $y_1$ in each $\akt{n, S} = \End(n)$, and a corresponding relation.)  Set $\kt{*, S, r} = \akt{*, S}/I_*$.  Then the algebra
 $\kt{*, S, r}$ is the algebra associated to a quotient category $KT(r)$, and 
 $\kt{n, S, r}$ (as previously defined)  can be identified with 
 $\End(n)$ in this category.

\subsection{Admissibility}  The cyclotomic BMW algebras and Kauffman tangle algebras can be defined over an arbitrary commutative unital ring $S$ with parameters $\rho$, $q$,  $\delta_j$    ($j \ge 0$), and
$u_1, \dots, u_r$ (roots of the polynomial satisfied by $y_1 = \rho x_1$), assuming   that $\rho$, $q$,  $\delta_0$  and  $u_1, \dots, u_r$   are invertible, 
and $\rho\inv - \rho=   (q\inv -q) (\delta_0 - 1)$.    However, unless the parameters
satisfy additional relations,  the identity element $\bm 1$ of the cyclotomic Kauffman tangle algebras
will be a torsion element over $S$;  if $S$ is a field (and the additional relations do not hold)  then $\bm 1 = 0$,  so $\kt {n, S, r}$ becomes trivial.    The additional conditions are called ``weak admissibility;"  see Section \ref{section: weak admissibility and the Markov trace} for details.

In order to obtain substantial results about the cyclotomic BMW and Kauffman tangle algebras, it seems necessary to impose stronger  conditions on the ground ring $S$.
An appropriate condition, known as ``admissibility,"  was introduced by Wilcox and Yu in ~\cite{Wilcox-Yu}.  
Their condition has a simple formulation in terms of the 
 representation theory of the 2--strand algebra $\bmw{2, S, r}$,  and also translates    into explicit relations on the parameters.  See Section \ref{section: admissiblity} for further details.

\subsection{Results} 
The main results of this paper and ~\cite{GH3} is that if the ground ring $S$ is an integral domain and admissible in the sense of Wilcox and Yu,  then $\bmw{n, S, r} \cong \kt{n, S, r}$,  and, moreover, these algebras are free $S$--modules of rank $ r^n (2n-1)!!$.  The proof of these results has a topological component,  which is given in this paper,  and an algebraic component, which is given in ~\cite{GH3}.

Our topological argument provides a straightening procedure for affine tangle diagrams that allows any affine tangle diagram to be expressed as a linear combination of affine tangle diagrams in a certain normal form (modulo Kauffman skein relations).

The straightening procedure has two important consequences: First,  it allows us to produce a spanning set $\A'_r$   of cardinality $r^n (2n-1)!!$ for the cyclotomic BMW algebra $\bmw{n, S, r}$. 
Second, it allows us to show that, if $S$  is weakly admissible,   then  $\kt{0, S, r}$ is a free $S$--module of rank 1.  Freeness of  $\kt{0, S, r}$  implies the existence of the ``Markov trace"  $\eps: \kt{n, S, r} \to \kt{0, S, r} \cong S,$
which is defined on the level of affine tangle diagrams by ``closing" diagrams:
$$
{\rm \eps} :  \inlinegraphic{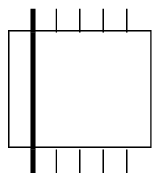} \quad \mapsto \quad  \delta_0^{-n}\ 
\inlinegraphic{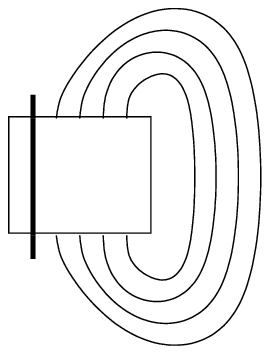}
$$
(Freeness of $\kt{0, S, r}$, or the existence of the Markov trace, immediately implies the existence of cyclotomic Kauffman link invariants in the solid torus,  cf.  ~\cite{Morton-Wassermann, LambrJKTR1999, Turaev-Kauffman-skein}.)  

To achieve the main results, it remains to show that the spanning set $\varphi(\A'_r)$
of the cyclotomic Kauffman tangle algebra $\kt{n, S, r}$ is linearly independent, when the ground ring $S$ is an admissible integral domain.   This is done in ~\cite{GH3}.  We give a brief summary of the strategy.  First, we show that there is a universal admissible integral domain $\overline R$, such that every admissible integral domain is a quotient of $\overline R$.  We then analyze the representation theory of the cyclotomic BMW algebras  defined over the field of fractions $F$ of $\overline R$, by adapting the inductive method of Wenzl from ~\cite{Birman-Wenzl, Wenzl-Brauer, Wenzl-BCD}.   This analysis shows that $\kt{n, F, r} \cong \bmw{n, F, r}$, and that the dimension of these algebras over $F$ is $r^n (2n-1)!!$.  A relatively simple argument (given in Section \ref{section: admissiblity} of this paper)  then shows that for any admissible integral domain $S$,  we have
$\kt{n, S, r} \cong \bmw{n, S, r}$, and these algebras are free of rank $r^n (2n-1)!!$ over $S$.

An outline of the paper is as follows:   In Section \ref{A  new basis of the affine Birman-Wenzl-Murakami algebra} we discuss the straightening procedure for affine tangle diagrams and its consequences.  In Section \ref{The cyclotomic BMW and Kauffman tangle algebras}, the cyclotomic algebras are introduced.  In Section \ref{section: weak admissibility and the Markov trace},   we discuss weak admissibility and the existence of the Markov trace.  In Section \ref{section: admissiblity}  we summarize ~\cite{GH3} and explain how our main theorem follows by combining the results of this paper with ~\cite{GH3}.
 In the final section of the paper,  we remark that a similar straightening procedure  can be applied to the affine and cyclotomic Hecke algebras,  allowing us to recover a result of Lambropoulou  ~\cite {LambrJKTR1999} with a less computationally intensive proof.

\subsection{Related work, and acknowledgments} Wilcox and Yu  have been studying the same material independently and have obtained  similar (and slightly stronger)  results ~\cite{Wilcox-Yu, Yu-thesis, Wilcox-Yu2}.   We are indebted to 
 Wilcox and Yu for pointing out an error in a previous preprint version of our work,  which required us to substantially rework the algebraic (linear independence) component of the arguments.  In fact, we had to adopt a completely different strategy for proving linear independence.  Meanwhile,  Wilcox and Yu ~\cite{Yu-thesis, Wilcox-Yu2} were able to make our original strategy work, using  a refined analysis of their admissibility condition. 

We would also like to mention here recent work of 
Ariki, Mathas and Rui on the 
``cyclotomic Nazarov--Wenzl algebras" ~\cite{ariki-mathas-rui},  which are cyclotomic quotients of the degenerate affine BMW algebras introduced by Nazarov ~\cite{Nazarov}.  Ideas from ~\cite{ariki-mathas-rui}  play an important role in the companion paper ~\cite{GH3}, and thus in our project as a whole.

It is shown in ~\cite{G-cellular}  and ~\cite{Yu-thesis}  that cyclotomic BMW algebras defined over integral admissible ground rings are cellular
in the sense of Graham and Lehrer ~\cite{Graham-Lehrer-cellular}.   The proof in ~\cite{G-cellular} relies on this paper and ~\cite{GH3}.

\section{A  new basis of the affine Birman-Wenzl-Murakami algebra}
\label{A  new basis of the affine Birman-Wenzl-Murakami algebra}

\subsection{The affine Kauffman tangle and BMW algebras}
We begin by recalling the definitions of the affine Birman-Wenzl-Murakami algebra and of the affine Kauffman tangle algebra.

The  {\em affine Kauffman tangle algebra}  $\akt n$ is the algebra of framed $(n, n)$--tangles in  $A \times I$, where $A$ is the annulus and $I$ the interval,  modulo the Kauffman skein relations.  This algebra can be described in terms of {\em affine tangle diagrams}, as follows.

\begin{figure} [ht]
\begin{eqnarray*}
\text{I}&\quad\ &\inlinegraphic{right_twist} \quad \longleftrightarrow
 \quad \inlinegraphic{vertical_line} \quad
  \longleftrightarrow  \quad \inlinegraphic{left_twist}\\
\text{II}&\quad\ &\inlinegraphic[scale =.5] {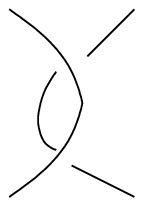} \quad 
\longleftrightarrow 
\quad \inlinegraphic[scale=1.75]{id_smoothing} \\
  \text{III}&\quad\ &\inlinegraphic{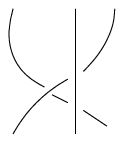} \quad \longleftrightarrow \quad 
\inlinegraphic{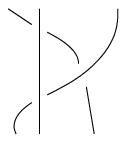} 
\end{eqnarray*} 
\caption{Reidemeister moves} \label{figure-Reidemeister moves}
\end{figure}

An ordinary   tangle diagram is a piece of knot diagram in the rectangle $\mathcal R = I \times I$  consisting of some number of closed curves and some number of topological intervals.   The intervals must have their endpoints on the upper or lower edge of the rectangle.  A $(k, n)$--tangle diagram is one with $k$ vertices (endpoints of intervals) on the upper edge of $\mathcal R$ and $n$ vertices on the lower edge.   We regard two ordinary tangle diagrams as equivalent if they are regularly isotopic (i.e. connected by a sequence of Reidemeister moves of types II and III,  followed by an isotopy of $\mathcal R$,  see Figure \ref{figure-Reidemeister moves}.   One can compose a $(k, n)$--tangle diagram $a$  and an  $(n, m)$--tangle diagram $b$  by ``stacking"  $a$ over $b$.  This yields a monoid structure on regular isotopy classes of $(n, n)$--tangle diagrams.

\begin{figure}[ht]
$$
\inlinegraphic[scale=1.5]{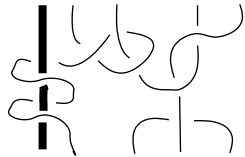}
$$
\caption{Affine tangle diagram}\label{figure-affine tangle diagram}
\end{figure}

An {\em affine  $(k,n)$--tangle diagram}\/ is an
ordinary
$(k+1, n+1)$--tangle diagram which includes a distinguished vertical curve.
We will draw affine tangle diagrams with the distinguished curve
drawn as a thickened vertical segment.   We refer to the distinguished
curve as the ``flagpole'',  and to the other curves in the diagram as ``ordinary strands". 
A picture of a typical affine tangle diagram is given in Figure \ref{figure-affine tangle diagram}.   Regular isotopy classes of affine $(n,n)$--tangle diagrams are closed under composition.  We let $\uhat n n $  denote the monoid of regular isotopy classes of affine $(n,n)$--tangle diagrams under composition.

For $j \ge 0$, let $\varTheta_j$  (resp. $\varTheta_{-j}$)  denote the (regular isotopy class of) the closed
curve with no self--crossings that winds $j$ times around the flagpole in the positive sense  (resp. in the negative sense). 
$$
\begin{array}{c c}
\includegraphics[scale=1.5]{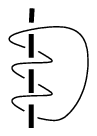} \\[-6pt] \vrule height20pt width0pt depth0pt 
{\scriptstyle \varTheta_3} \end{array} 
\qquad
\begin{array}{c}
\includegraphics[scale=1.5]{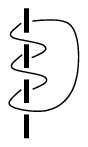}  \\[-6pt]  \vrule height18pt width0pt depth0pt  {\scriptstyle \varTheta_{-3}} 
\end{array}
$$
Note that $\varTheta_0$ is represented by a closed curve that does not intersect the flagpole.

\begin{definition}\rm\label{definition affine Kauffman tangle algebra}
Let $S$ be a commutative unital ring containing   elements 
$\rho$, $q$, and  $\delta_j$, $j \ge 0$,   with $\rho$, $q$,  and $\delta_0$ invertible, satisfying the relation
$
\rho\inv - \rho=  \break (q\inv -q) (\delta_0 - 1).
$
The {\em affine Kauffman tangle
algebra} 
$
\akt {n, S} = \akt {n, S}(\rho, q, \delta_0, \delta_1, \dots)
$
is the monoid algebra $S\,  \uhat n n$ modulo the following relations:
\begin{enumerate}
\item (Crossing relation)
$$
\qquad\quad \inlinegraphic[scale=.6]{pos_crossing} - \inlinegraphic[scale=.3]{neg_crossing} 
\quad = 
\quad
(q\inv - q)\,\left( \inlinegraphic[scale=1.5]{e_smoothing} - 
\inlinegraphic[scale=1.5]{id_smoothing}\right).
$$
\item (Untwisting relation)
$$\qquad\quad 
\inlinegraphic{right_twist} \quad = \quad \rho \quad
\inlinegraphic{vertical_line} \quad\ \text{and} \quad\ 
\inlinegraphic{left_twist} \quad = \quad \rho\inv \quad
\inlinegraphic{vertical_line}. 
$$
\item  (Free loop relations ) For $j \ge 0$,  $T\, \cup \, \varTheta_j =  \rho^{-j}\delta_j T,$
 where $T\, \cup\,\ \varTheta_j$ is the union
of an affine  tangle diagram $T$ and a copy of the curve $\varTheta_j$, such that  there are no crossings between $T$ 
and~$\varTheta_j$.
\end{enumerate}
\end{definition}

\begin{remark}
The idea behind relation (3) is the following:  If one only imposes  relations (1) and (2),  
then the  $(0,0)$--affine tangle algebra is the polynomial algebra generated by  $\varTheta_j$   ($j \ge 0$),  over whatever ground ring one is working, and embeds in the center of the $(n,n)$--tangle algebra.   Therefore,  it makes sense to absorb the 
 $(0,0)$--affine tangle algebra into the ground ring, and this is accomplished by relation (3).
 \end{remark}

We now introduce the affine Birman--Wenzl--Murakami  (BMW) algebra. As above, let 
 $S$ be a commutative unital ring with
  elements 
$\rho$, $q$, and  $\delta_j$, $j \ge 0$,   with $\rho$, $q$,  and $\delta_0$ invertible, satisfying the relation
$
\rho\inv - \rho=   (q\inv-q) (\delta_0 - 1).
$

\begin{definition}\rm\label{definition affine BMW}
 The {\em affine
Birman--Wenzl--Murakami} algebra
$\abmw  {n, S}$ is the
$S$ algebra with generators $y_1^{\pm 1}$, $g_i^{\pm 1}$  and
$e_i$ ($1 \le i \le n-1$) and relations:
\begin{enumerate}
\item (Inverses)\hods $g_i g_i\inv = g_i\inv g_i = 1$ and 
$y_1 y_1\inv = y_1\inv y_1= 1$.
\item (Idempotent relation)\hods $e_i^2 = \delta_0 e_i$.
\item (Type $B$ braid relations) 
\begin{enumerate}
\item[\rm(a)] $g_i g_{i+1} g_i = g_{i+1} g_ig_{i+1}$ and 
$g_i g_j = g_j g_i$ if $|i-j|  \ge 2$.
\item[\rm(b)] $y_1 g_1 y_1 g_1 = g_1 y_1 g_1 y_1$ and $y_1 g_j =
g_j y_1 $ if $j \ge 2$.
\end{enumerate}
\item[\rm(4)] (Commutation relations) 
\begin{enumerate}
\item[\rm(a)] $g_i e_j = e_j g_i$  and
$e_i e_j = e_j e_i$  if $|i-
j|
\ge 2$. 
\item[\rm(b)] $y_1 e_j = e_j y_1$ if $j \ge 2$.
\end{enumerate}
\item[\rm(5)] (Affine tangle relations)\vadjust{\vskip-2pt\vskip0pt}
\begin{enumerate}
\item[\rm(a)] $e_i e_{i\pm 1} e_i = e_i$,
\item[\rm(b)] $g_i g_{i\pm 1} e_i = e_{i\pm 1} e_i$ and
$ e_i  g_{i\pm 1} g_i=   e_ie_{i\pm 1}$.
\item[\rm(c)\hskip1.2pt] For $j \ge 1$, $e_1 y_1^{ j} e_1 = \delta_j e_1$. 
\vadjust{\vskip-
2pt\vskip0pt}
\end{enumerate}
\item[\rm(6)] (Kauffman skein relation)\hods  $g_i - g_i\inv = (q\inv-q)(e_i -1)$.
\item[\rm(7)] (Untwisting relations)\hods $g_i e_i = e_i g_i = \rho \inv e_i$
 and $e_i g_{i \pm 1} e_i = \rho  e_i$.
\item[\rm(8)] (Unwrapping relation)\hods $e_1 y_1 g_1 y_1 = \rho e_1 = y_1 
g_1 y_1 e_1$.
\end{enumerate}
\end{definition}

\begin{remark}  The presentation differs  slightlyÊ\ from the one we used in ~\cite{GH1}.    There we used the generator $x_1 = \rho \inv y_1$,   and
 parameters $\vartheta_j =  \rho^{-j} \delta_j$    (so that $e_1 x_1^j e_1 = \vartheta_j e_1$).   We also used the parameter $z$  in place of $(q\inv -q)$.
  \end{remark}
  
Let $X_1$,  $G_i$,  $E_i$  denote the following affine tangle diagrams:
$$
X_1 = \inlinegraphic[scale= .7]{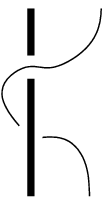}
\qquad
G_i =  \inlinegraphic[scale=.6]{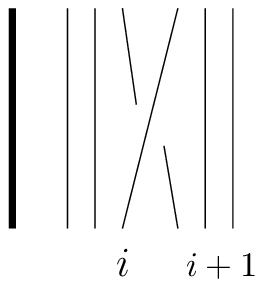}\qquad
E_i =  \inlinegraphic[scale= .7]{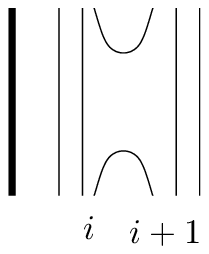} 
$$

\begin{theorem} [\cite{GH1}]
\label{theorem:  isomorphism affine BMW and KT}
 Let $S$ be any commutative unital ring with distinguished elements
$\rho$,  $q$,  and $\delta_j$,  $j \ge 0$,  as above.  
The affine BMW algebra $\abmw {n, S}$ is isomorphic to the affine Kauffman tangle algebra $\akt {n, S} $ by a map $\varphi$ determined by
$\varphi(g_i) = G_i$,  $\varphi(e_i) = E_i$,   and $\varphi(y_1) =   \rho X_1$.
\end{theorem}

\subsection{Relations in the $(0,0)$--tangle algebra}

\begin{figure}[t] 
$$
\begin{array}{c c}
\includegraphics[scale=1.3]{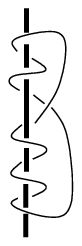} \\[-6pt] \vrule height20pt width0pt depth0pt 
{\scriptstyle \varTheta_{2, 3}} \end{array} 
\qquad
\begin{array}{c}
\includegraphics[scale=1.3]{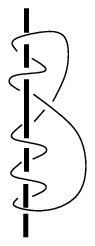}  \\[-6pt]  \vrule height18pt width0pt depth0pt  {\scriptstyle \varTheta_{2, 3}^{-}} 
\end{array}
$$
\caption{}   \label{figure: Theta a b}
\end{figure}

Since the  $(0,0)$--affine tangle algebra is generated by  $\varTheta_j$,  $j \ge 0$, 
in particular each 
$\varTheta_{-k} $ can be expressed as a polynomial in $\varTheta_{j}$, $j \ge 0$.   
We will find a recursive formula for $\varTheta_{-k}$   (correcting  a minor error in ~\cite{GH1},  Lemma 2.4.)
For $a \ge 1, b \ge 0$,  let $\varTheta_{a, b}$ be the curve with $a$ positive 
windings around the flagpole, and one positive crossing, and $b$ negative windings;   let  $\varTheta_{a, b}^{-}$ be the curve with the crossing reversed;  see Figure \ref{figure: Theta a b}.

\begin{lemma}  \label{lemma - recursion for f_r} 
Let $j \ge 1$.
\begin{enumerate}
\item $\varTheta_{-j} = \rho\, \varTheta_{1, j-1}.$
\item $\varTheta_j = \rho\,\inv\varTheta_{j-1, 1}^{-}.$
\item  For $a \ge 1, b\ge 2$,  $\varTheta_{a, b} =  \rho^2 \varTheta_{a+1, b-1} + (q\inv -q)(\varTheta_{a}\varTheta_{-b} - \varTheta_{a-b}).
$
\item  
$
\varTheta_{-j} = \displaystyle \rho^{2j-2}  \varTheta_j   + (q\inv - q) \sum_{k = 1}^{j-1} \rho^{2k -1}  (\varTheta_k \varTheta_{k-j} - \varTheta_{2k - j}).
$
\end{enumerate}
\end{lemma}

\begin{proof}  Point (1) follows from introducing
 a twist at the top of $\varTheta_{-j}$:
 $$
\inlinegraphic[scale=1.5]{Theta3inv} = \rho \inlinegraphic[scale=1.5]{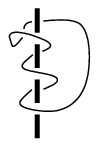} = \rho \inlinegraphic[scale=1.5]{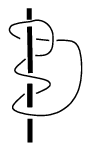}.
 $$
 The proof of (2) is the same.
 The tangle  obtained by smoothing the crossing  in $\varTheta_{a, b}$
horizontally is  $\varTheta_a \varTheta_{-b}$,  and the tangle obtained by smoothing the crossing vertically is
$\varTheta_{a-b}$, while $\varTheta_{a, b}^{-} = \rho^2 \varTheta_{a+1, b-1}$.
Thus the Kauffman skein relation gives
$
\varTheta_{a, b} =  \rho^2 \varTheta_{a+1, b-1} + (q\inv - q)(\varTheta_{a}\varTheta_{-b} - \varTheta_{a-b}).
$
An induction based on points (1)--(3) yields (4).
\end{proof}

\subsection{Flagpole--descending affine tangle diagrams}

\begin{definition}\rm
A {\em simple winding} is a piece of an affine tangle
diagram with one ordinary strand, without self--crossings, 
 regularly isotopic to the intersection of one of the
affine tangle diagrams
$X_1$ or $X_1\inv$ with a small neighborhood of the flagpole, as in  the following
figure.
\end{definition}

\centerline{$\inlinegraphic{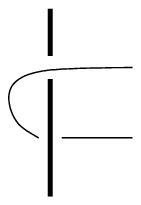}$ }

\begin{figure}[ht] 
\centerline{$\inlinegraphic{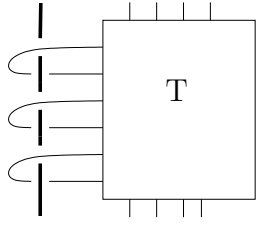}$}
\caption{Affine tangle diagram in standard position} \label{figure-standard position}
\end{figure}

\begin{definition}\rm
\label{definition: standard position}
An affine tangle diagram is in {\em standard
position}  if:
\begin{enumerate}
\item  It has no crossings to the left of the flagpole.
\item  There is a neighborhood of the flagpole whose intersection with
the tangle diagram is a union of simple windings.
\item  The simple windings have no crossings and are not nested.  That is,
between the two crossings of a simple winding with the flagpole, there is
no other crossing of a strand with the flagpole.
\end{enumerate}
\end{definition}
See Figure \ref{figure-standard position}.

\begin{lemma}[\cite{GH1}]
Any affine tangle diagram is regularly isotopic to 
an affine tangle diagram in standard position.
\end{lemma}

Fix an integer $n \ge 1$.
Recall that an  $(n,n)$--affine tangle diagram is  a  figure contained in $ I \times I$.
For fixed points $0 < a_1 < \cdots < a_n < 1$ in $I$,  we write
$\p i = (a_i, 1)$ and $\overline{ \p i} = (a_i, 0)$;  we refer to these points as {\em vertices}.   We order the vertices   of $(n,n)$-- affine tangle diagrams by
$
\p 1 < \p 2 < \cdots < \p n < \pbar n  < \cdots < \pbar 2 < \pbar 1.
$
That is, the vertices are taken in clockwise order around the boundary of $I \times I$.

\begin{definition}\rm
An {\em orientation} of an  affine or ordinary 
  tangle diagram  is
  a linear ordering of the strands,
 a choice of an orientation of each strand, and a choice of an initial point on 
each closed loop.
   \end{definition}
 
 An orientation determines a way of traversing the tangle diagram;  namely,  
the 
strands are traversed successively, in
the given order and orientation (the closed loops being traversed starting at 
the assigned initial point).

\begin{definition}\rm
A {\em standard orientation} of an ordinary or affine 
$(n,n)$--tangle diagram is one in
which 
\begin{enumerate}
\item each non-closed strand is oriented from its lower numbered vertex to its 
higher numbered vertex.
\item    The non-closed strands precede the closed loops.
\item  The  non-closed strands  are ordered according to the order of the initial vertices.
\end{enumerate}
\end{definition}

If a tangle diagram has no closed loops, then it has a unique standard 
orientation.

  \begin{definition}\rm
  An oriented affine or ordinary tangle diagram is {\em stratified} if 
  \begin{enumerate}
  \item
  there is a linear ordering of the strands  such that if strand $s$ precedes strand $t$ in the order, then each crossing of $s$ with $t$ is an over--crossing. 
   \item  each strand is totally descending,  that is, each self--crossing of the strand is encountered first as an over--crossing as the strand is traversed according to the orientation.
   \end{enumerate}  
   We call the corresponding ordering of the strands the {\em stratification order}.  
  \end{definition}
  
  Note that a stratification order need not coincide with the ordering of strands determined by the orientation.

  \begin{lemma}\label{lemma-totally descending tangles
span}
Endow affine $(n,n)$--tangles with an orientation and a stratification order, each determined by some rule depending  on the vertices of the strands.
$\akt {n, S}$ is spanned by affine $(n,n)$--tangle diagrams in standard position
 that are  stratified with respect to the given  stratification order.
  \end{lemma}
  
  \begin{proof}    We already know that $\akt{n, S}$ is spanned by affine tangle
  diagrams in standard position.   Observe that changing or smoothing a crossing
  of ordinary strands leaves such a diagram in standard position.  In the following
  paragraph ``tangle" means  ``affine tangle diagram in standard position".
  
 The proof is by induction on the number of crossings of ordinary strands.
  If a tangle has no crossings, it is already stratified.
  Let $T$ be a tangle  with $ l  \ge 1$ crossings, and  assume that any 
tangle  with fewer
than $ l $ crossings is in the span of stratified tangles.   
The stratified
tangle $S$ which differs from $T$ only by reversing some number of 
crossings
  is congruent to $T$ modulo the span of tangles with fewer crossings, 
hence modulo the span of
stratified tangles.
  \end{proof}

  \begin{definition}  An oriented, stratified affine tangle diagram $T$  in standard position is said to be {\em flagpole descending} if it satisfies the following conditions:
  \begin{enumerate}
     \item  $T$ is not regularly isotopic to an affine tangle diagram in standard position with fewer simple windings.
  \item  The strands of $T$ have no self--crossings.
  \item   As  $T$ is traversed according to the orientation,  successive
  crossings of ordinary strands with the flagpole descend the flagpole.   
  \end{enumerate}
  \end{definition}

\begin{figure}[ht]
$$
\begin{aligned}[c]\setlength{\unitlength}{0.00062500in}
\begingroup\makeatletter\ifx\SetFigFont\undefined%
\gdef\SetFigFont#1#2#3#4#5{%
  \reset@font\fontsize{#1}{#2pt}%
  \fontfamily{#3}\fontseries{#4}\fontshape{#5}%
  \selectfont}%
\fi\endgroup%
{\renewcommand{\dashlinestretch}{30}
\begin{picture}(815,1281)(0,-10)
\thicklines
\path(278,783)(278,33)
\path(278,1233)(278,933)
\thinlines
\path(353,408)(803,408)
\path(803,858)(800,858)(794,858)
	(782,858)(765,858)(741,857)
	(712,857)(678,857)(641,856)
	(602,855)(563,855)(524,854)
	(486,853)(451,852)(417,850)
	(387,849)(358,848)(333,846)
	(309,844)(287,843)(267,841)
	(249,838)(232,836)(216,833)
	(193,828)(171,823)(151,817)
	(133,811)(115,804)(99,796)
	(84,787)(70,778)(58,768)
	(48,758)(39,748)(32,738)
	(26,727)(22,717)(18,706)
	(16,695)(13,681)(12,666)
	(13,649)(15,632)(18,615)
	(22,597)(28,580)(34,563)
	(41,548)(49,533)(57,520)
	(66,508)(74,497)(84,487)
	(95,477)(108,467)(123,456)
	(140,445)(158,434)(175,424)(203,408)
\path(122.998,419.163)(203.000,408.000)(152.766,471.258)
\end{picture}
} \\ {\rm(a)} \\ \end{aligned} 
\qquad \qquad
\begin{aligned}[c] \setlength{\unitlength}{0.00062500in}
\begingroup\makeatletter\ifx\SetFigFont\undefined%
\gdef\SetFigFont#1#2#3#4#5{%
  \reset@font\fontsize{#1}{#2pt}%
  \fontfamily{#3}\fontseries{#4}\fontshape{#5}%
  \selectfont}%
\fi\endgroup%
{\renewcommand{\dashlinestretch}{30}
\begin{picture}(815,1281)(0,-10)
\thicklines
\path(278,783)(278,33)
\path(278,1233)(278,933)
\thinlines
\path(803,858)(800,858)(794,858)
	(782,858)(765,858)(741,857)
	(712,857)(678,857)(641,856)
	(602,855)(563,855)(524,854)
	(486,853)(451,852)(417,850)
	(387,849)(358,848)(333,846)
	(309,844)(287,843)(267,841)
	(249,838)(232,836)(216,833)
	(193,828)(171,823)(151,817)
	(133,811)(115,804)(99,796)
	(84,787)(70,778)(58,768)
	(48,758)(39,748)(32,738)
	(26,727)(22,717)(18,706)
	(16,695)(13,681)(12,666)
	(13,649)(15,632)(18,615)
	(22,597)(28,580)(34,563)
	(41,548)(49,533)(57,520)
	(66,508)(74,497)(84,487)
	(95,477)(108,467)(123,456)
	(140,445)(158,434)(175,424)
	(189,416)(198,411)(202,408)(203,408)
\path(473.000,438.000)(353.000,408.000)(473.000,378.000)
\path(353,408)(803,408)
\end{picture}
} \\ {\rm(b)} \\ 
\end{aligned}
\qquad \qquad
\begin{aligned}[c]\setlength{\unitlength}{0.00062500in}
\begingroup\makeatletter\ifx\SetFigFont\undefined%
\gdef\SetFigFont#1#2#3#4#5{%
  \reset@font\fontsize{#1}{#2pt}%
  \fontfamily{#3}\fontseries{#4}\fontshape{#5}%
  \selectfont}%
\fi\endgroup%
{\renewcommand{\dashlinestretch}{30}
\begin{picture}(692,1281)(0,-10)
\thicklines
\path(230,258)(230,33)
\path(230,1233)(230,408)
\thinlines
\path(155,783)(154,783)(150,780)
	(141,775)(128,766)(111,756)
	(94,744)(78,732)(64,720)
	(53,708)(44,697)(36,684)
	(30,670)(25,658)(21,644)
	(18,628)(15,612)(14,595)
	(12,577)(12,558)(12,539)
	(14,521)(15,504)(18,488)
	(21,472)(25,458)(30,445)
	(35,435)(40,425)(46,416)
	(53,407)(61,399)(70,391)
	(81,384)(92,377)(105,371)
	(118,366)(133,361)(148,357)
	(164,354)(181,350)(199,348)
	(218,345)(233,344)(250,342)
	(269,341)(289,340)(311,339)
	(336,338)(364,337)(394,336)
	(427,336)(462,335)(499,335)
	(536,334)(571,334)(604,334)
	(632,333)(653,333)(668,333)
	(676,333)(679,333)(680,333)
\path(680,783)(305,783)
\path(395.000,813.000)(305.000,783.000)(395.000,753.000)
\end{picture}
} \\ {\rm(c)} \\ 
\end{aligned} 
\qquad \qquad
\begin{aligned}[c]\setlength{\unitlength}{0.00062500in}
\begingroup\makeatletter\ifx\SetFigFont\undefined%
\gdef\SetFigFont#1#2#3#4#5{%
  \reset@font\fontsize{#1}{#2pt}%
  \fontfamily{#3}\fontseries{#4}\fontshape{#5}%
  \selectfont}%
\fi\endgroup%
{\renewcommand{\dashlinestretch}{30}
\begin{picture}(692,1281)(0,-10)
\thicklines
\path(230,258)(230,33)
\path(230,1233)(230,408)
\thinlines
\path(680,783)(305,783)
\path(69.436,693.675)(155.000,783.000)(37.468,744.449)
\path(155,783)(128,766)(111,756)
	(94,744)(78,732)(64,720)
	(53,708)(44,697)(36,684)
	(30,670)(25,658)(21,644)
	(18,628)(15,612)(14,595)
	(12,577)(12,558)(12,539)
	(14,521)(15,504)(18,488)
	(21,472)(25,458)(30,445)
	(35,435)(40,425)(46,416)
	(53,407)(61,399)(70,391)
	(81,384)(92,377)(105,371)
	(118,366)(133,361)(148,357)
	(164,354)(181,350)(199,348)
	(218,345)(233,344)(250,342)
	(269,341)(289,340)(311,339)
	(336,338)(364,337)(394,336)
	(427,336)(462,335)(499,335)
	(536,334)(571,334)(604,334)
	(632,333)(653,333)(668,333)
	(676,333)(679,333)(680,333)
\end{picture}
} \\ {\rm(d)} \\ 
\end{aligned}.
$$
\caption{Types of oriented simple windings}  \label{figure-simple winding types}
\end{figure}
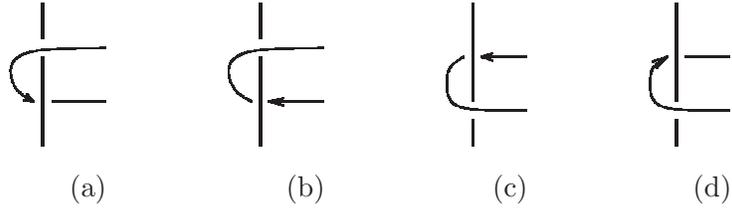

  \begin{remark}  \label{remark:  on flagpole descending diagrams and simple windings}
   In a flagpole descending affine tangle diagram:
  \begin{enumerate}
  \item Of the four types of oriented simple windings (see Figure \ref{figure-simple winding types}), only types
  (a) and (c)  can occur.   Moreover,  a strand can have windings of one type only,  since otherwise the number of simple windings could be reduced by regular isotopy.
\item  For each strand $s$ there is a neighborhood  $N$  of an interval on the flagpole such that:   $N$ does not intersect any strand other than $s$,   all crossings of $s$ with the flagpole are contained in $N$,  and
  the intersection of $s$ with $N$  is regularly isotopic to  $X_1^k $ for some $k \in \Z$.
\end{enumerate}
  \end{remark}
  
Our immediate goal is to show that $\akt{n, S}$ is spanned by affine tangle diagrams without closed strands that are stratified (with an arbitrary stratification order) and flagpole descending with respect to the standard orientation.

\begin{figure}[ht]  
\centerline{$\inlinegraphic{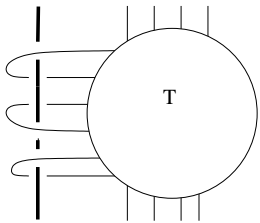}$}
\caption{Affine tangle diagram in circular form} \label{figure-circular form}
\end{figure}

 Let $D$ be the circle of radius $1/3$, centered at  $((a_1 + a_n)/2, 1/2)$.   Planar isotopy alone can transform any affine $(n,n)$--tangle diagram $T$ in standard position so that 
 \begin{itemize}
 \item $D$ lies transversal to $T$, 
 \item each vertex is connected by a straight line segment to a point of $D$, 
 \item no crossings of ordinary strands occur outside of $D$, and
 \item   the part of the affine tangle diagram outside of $D$ consists only of simple windings beginning and ending on  $D$ and line segments connecting  $D$ with the vertices.  
\end{itemize}
 See Figure \ref{figure-circular form}.  Call a representative affine tangle diagram in this form a {\em circular form} of the affine tangle diagram  $T$.  

Let $T$ be an  $(n,n)$--affine tangle diagram in circular form.  If $k$ is the number of
simple windings in $T$,  then  $T$ intersects $D$ in $2(n+k)$ points.
  The  $2k$ points of intersection incident with  simple windings will be called {\em winding points} and those connected by a line segment  with a vertex will be called {\em vertex points}.   
The   $2(n+k)$ points of interesection of $D$ with $T$  are connected in pairs by
curves running through the interior of $D$,  which we call $D$--arcs.

\begin{definition} A {\em tight} affine tangle diagram is an affine tangle diagram in circular form  is which every $D$--arc is  straight line segment.
\end{definition}

\begin{lemma}  Every stratified affine tangle diagram $T$  in standard position is ambient isotopic to a tight affine tangle diagram  $T'$.    Moreover,  $T'$ has at most as many crossings as $T$.
\end{lemma}
 
\begin{proof}  As observed above, a stratified affine tangle diagram in standard position is equivalent by planar isotopy to a diagram $T$  in circular form.   Since $T$
is stratified,  its $D$--arcs can be lifted slightly out of the plane so that they lie
at different levels above the plane.  Each $D$--arc is unknotted, so it can be changed by ambient isotopy to a straight line segment.   (It may be necessary to adjust the positions of intersections of $T$ with $D$ to avoid triple intersections of $D$--arcs.)

Let $T$ and $T'$ be ambient isotopic affine tangle diagrams in circular form, with
$T'$ tight.    Let $s, t$ be two $D$--arcs in $T$  and let $s', t'$ be the corresponding $D$--arcs in $T'$  (with the same endpoints on $D$).      Then  $s'$ and $t'$  have exactly one crossing in $T'$  if the endpoints of $t'$ lie on opposite sides of $s'$  (and no crossings otherwise), while  $s$ and $t$ have at least one crossing if the endpoints of
$t$  lie on opposite sides of $s$.     Thus $T'$   has at most as many crossings as $T$.
\end{proof}

Consider an affine $(n,n)$--tangle diagram $T$ in circular form,  with $k$ simple windings, stratified with respect to  a standard orientation and some stratification order.
The $2n$  vertex points and the $2k$ winding points of $T$ lie on disjoint arcs of the circle $D$.    Let  $p_1, \dots, p_{2k}$  be the  winding points listed according to their position on $D$,  in counterclockwise order.     Let $w_1, \dots, w_{2k}$  be the winding points listed according to their order in the orientation of $T$.  There is a permutation
$\sigma = \sigma(T)$ of $\{1, \dots, 2k\}$ such that $w_i  = p_{\sigma(i)}$ for all $i$.    For 
$T$ to be flagpole descending,  it is necessary that $\sigma$ be the identity permutation;  conversely,
if $\sigma$ is the identity permutation, then $T$ can be reduced to a flagpole descending tangle diagram by ambient isotopy.   Moreover,  closed loops can be eliminated from a 
flagpole descending tangle by use of the free loop relation (3) in Definition  \ref{definition affine Kauffman tangle algebra}   (and Lemma \ref{lemma - recursion for f_r}.)

\begin{proposition}\label{proposition:  layered diagrams span}  Endow affine $(n,n)$--tangle diagrams with  any stratification order.  
$\akt{n, S}$ is spanned by  affine tangle diagrams  without closed loops that are
stratified according to the given stratification order and
 flagpole descending with respect to the standard orientation.
\end{proposition}

\begin{proof}  We know that $\akt{n, S}$ is spanned by affine tangle diagrams in standard position that are stratified  with respect  the given stratification order.

Let $T$ be a stratified affine tangle diagram in standard position with $m$ crossings of ordinary strands.
Let $k$ be the number of simple windings of $T$, and let  $p_i$,  $w_i$ \break ($1 \le i \le 2k$),
and $ \sigma(T)$ be as in the discussion preceding the statement of the Proposition.  
Write $\ell(\sigma)$  for the length of a permutation $\sigma$.

If $\sigma(T)$ is the identity permutation,  then  $T$ can be reduced to 
a flagpole descending affine tangle  diagram without closed strands,  by ambient isotopy and the free loop relations.

We claim that if   $\sigma(T)$ is not the identity permutation,  then
$T$ can be written as a sum $T = T' + T''$,  where $T'$ is an affine tangle diagram with  no more than $m$ crossings and $\ell(\sigma(T')) < \ell(\sigma(T))$;  and $T''$ is a linear combination of affine tangle diagrams with strictly fewer than $ m$ crossings.

If the claim is established,  then the result follows by a double induction.  Namely,  if $m = 0$,  then by the claim, $T$ is equal to an affine tangle  diagram $T'$ with no crossings and with $\ell(\sigma(T')) < \ell(\sigma(T))$.   It follows by induction on the length of $\sigma(T)$ that $T$ is equal to a multiple of a stratified, flagpole descending  affine tangle diagram with   no closed loops.
Now suppose that $T$ has $m \ge 1$ crossings and that  $\ell(\sigma(T)) > 1$;  assume inductively  that any stratified  affine tangle diagram $S$ with strictly fewer than $m$ crossings and also any stratified  affine tangle diagram  $S$ with $m$ crossings but with
$\ell(\sigma(S)) < \ell(\sigma(T))$ is in the span of stratified flagpole descending affine tangle diagrams without closed loops.    Now the assertion of the proposition follows immediately from this induction hypothesis and the claim. 

We proceed with the proof of the claim.  Since $\sigma(T)$ is not the identity permutation, let 
 $i_0$ be the first index such that $\sigma(i_0) > i_0$.   
 We consider two cases:

 \begin{figure}[h]
$$
\inlinegraphic[scale=.75]{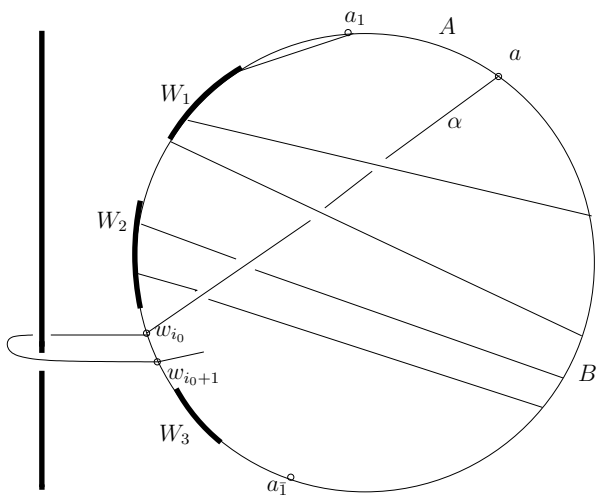}
$$
\caption{Sliding a winding, Case 1 before.}\label{figure slide before}
\end{figure}

{\em Case 1. The winding point $w_{i_0}$ is the first winding point on a non--closed strand $s$.}   We will suppose  that the oriented winding associated with the winding point $w_{i_0}$ is of type (c)  from Figure \ref{figure-simple winding types};   the other types are handled similarly.    Let $W_1$ be the set of winding points $w_i$ with $i < i_0$.
Let $W_2$ be the set of winding points  $w_i$  with $i \ge i_0 +2$ and $\sigma(i) < \sigma(i_0)$.    Let $W_3$  be the set of winding points $w_i$  with 
 $i \ge i_0 +2$ and $\sigma(i) > \sigma(i_0+1)$.    See Figure \ref{figure slide before}.
 Since $w_{i_0}$ is the first winding point on $s$,  it is connected by a $D$--arc $\alpha$ to a vertex point $a$.

 We order the vertex points  according to their position on $D$,  in clockwise order.  (The vertex points are in the same order as their associated vertices.)  Let $A$  be the set of vertex points that precede $a$  and let $B$  be the set of vertex points that follow $a$ in this order.

 The points of $A \cup W_1$  can be joined by $D$--arcs to points of $A \cup W_1$  or $B$.
 $D$--arcs that join $A\cup W_1$ and $B$  cross  $\alpha$.      The points of 
 $W_2$   can be joined by $D$--arcs to points of $W_2$,  or  $ \{w_{i_0 +1}\}  \cup W_3\cup B$.  $D$--arcs that join $W_2$ with $ \{w_{i_0 +1}\}  \cup W_3\cup B$ cross  $\alpha$.  This accounts for all $D$--arcs that cross $\alpha$.
 
  \begin{figure}[ht]
$$
\inlinegraphic[scale=.75]{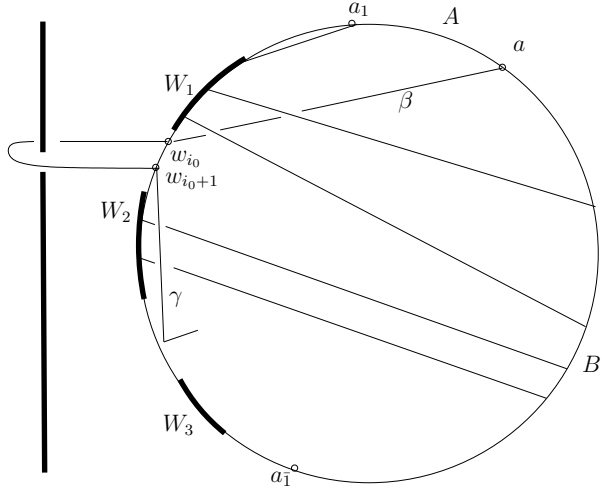}
$$
\caption{Sliding a winding, Case 1 after.}\label{figure slide after}
\end{figure}

 Let $T^{(1)}$ be the tight tangle diagram  obtained from $T$ by changing  all crossings of $\alpha$ with other $D$--arcs incident with $W_2$  to under--crossings.  $T$ and $T^{(1)}$  are congruent modulo the span of diagrams with fewer crossings.  $T^{(1)}$  is regularly isotopic to the diagram $T^{(2)}$ in which the
 simple winding associated with $w_{i_0}$ is slid to a position between
 $W_1$ and $W_2$,  as shown in Figure \ref{figure slide after}.   In this figure,
 $\beta$ is a $D$--arc connecting $a$ with (the new position of)  $w_{i_0}$,  and
 $\gamma$ is a curve connecting (the new position of)  $w_{i_0 +1 }$  with a point
 $x$ in the interior of $D$ close to the old position of $w_{i_0 +1 }$.  Finally let 
 $T^{(3)}$ be the tangle diagram obtained from  $T^{(2)}$ by changing crossings of 
 $\beta \cup \gamma$ with other strands to agree with the stratification order.  Then
  $T^{(2)}$ and  $T^{(3)}$ are congruent modulo the span of affine tangle diagrams with fewer crossings, as above.

 The $D$--arcs that cross the $D$--arc  $\beta$ in  $T^{(3)}$  are incident with $W_1 \cup A$;  they correspond
 one--to--one with $D$--arcs in $T$  incident with  $W_1 \cup A$ that cross $\alpha$.
The $D$--arcs that cross the curve $\gamma$ in  $T^{(3)}$  are incident with $W_2$;  they correspond
 one--to--one with $D$--arcs in $T$  incident with  $W_2$ that cross $\alpha$, with the following exception:   $T$ might have a $D$--arc connecting $w_{i_0 +1}$ with a point of 
 $W_2$, and crossing $\alpha$;   in $T^{(3)}$,  this $D$--arc is replaced by a curve
 connecting the endpoint $x$ of $\gamma$ with the point of $W_2$. 
 This discussion shows that $T^{(3)}$ has at most as many crossings as $T$.
Moreover,   $\ell(\sigma(T^{(3)})) < \ell(\sigma(T))$.

 {\em Case 2.  The winding point $w_{i_0}$  is on a closed loop;  or $w_{i_0}$ is on a non--closed strand $s$, but it  is not the first winding point on  $s$.}  
 
Recall that the order of closed loops and the orientation and initial point on each closed loop can be chosen arbitrarily in a standard orientation.   By making  appropriate choices, we can assume the following without loss of generality:

\begin{enumerate}  
\item    Suppose that $o$ and $o'$ are two closed loops in $T$ and $o$ precedes $o'$ in the ordering of strands.   If $w_i$  and $w_j$ are the  first  winding points on $o$ and $o'$, respectively,  with respect to the orientation,  then  $\sigma(i)  < \sigma(j)$.
\item   If $o$ is a closed loop in $T$,  $w_i$  is the first winding point  on $o$ and $w_j$ is another winding point on $o$,   then $\sigma(i)  < \sigma(j)$.
\end{enumerate} 
 
With these assumptions, if $w_{i_0}$  is on a closed loop $s$, then  $w_{i_0}$  is not the first winding point on $s$.

 Again, we suppose  that the oriented winding associated with the winding point $w_{i_0}$ is of type (c)  from Figure \ref{figure-simple winding types}, the   other types being handled similarly.  Define $W_1$,  $W_2$,  and $W_3$ as in Case 1(a).
 Since $w_{i_0}$ is not the first winding point on its strand,  it is connected by a $D$--arc $\alpha$ to  the last winding point in $W_1$.  All   $D$--arcs that cross $\alpha$ are incident with $W_2$, 
  see Figure \ref{figure case1b before}.
\begin{figure}[ht]
$$
\inlinegraphic[scale=.75]{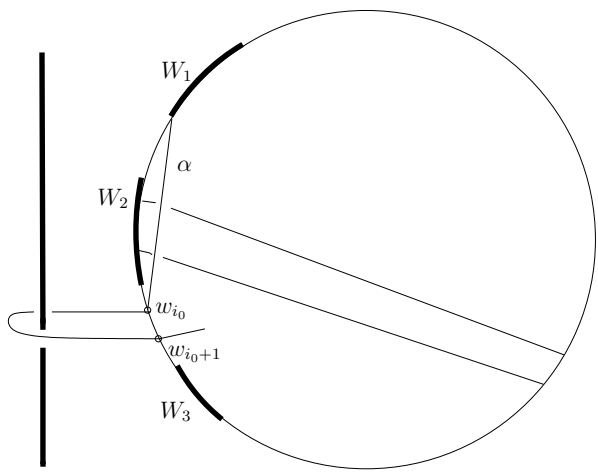}
$$
\caption{Sliding a winding, Case 2, before.}\label{figure case1b before}
\end{figure}

Now we  change all the crossings of $\alpha$ with other $D$--arcs to under--crossings,   and slide the simple winding associated with $w_{i_0}$ 
  to a position between  $W_1$ and $W_2$, to obtain a diagram  $T^{(2)}$,   as shown in Figure \ref{figure case1b after}.    Let $T^{(3)}$  be the tangle diagram obtained from 
  $T^{(2)}$ by changing crossings of $\alpha$ with other strands to agree with the stratification order.

  Then $T$ is regularly isotopic with $ T^{(3)}$ modulo
  the span of diagrams with fewer crossings.  Moreover,  $ T^{(3)}$ has at most as many crossings as $T$, 
 and  $\ell(\sigma( T^{(3)})) < \ell(\sigma(T))$.   
   \begin{figure}[ht]
$$
\inlinegraphic[scale=.75]{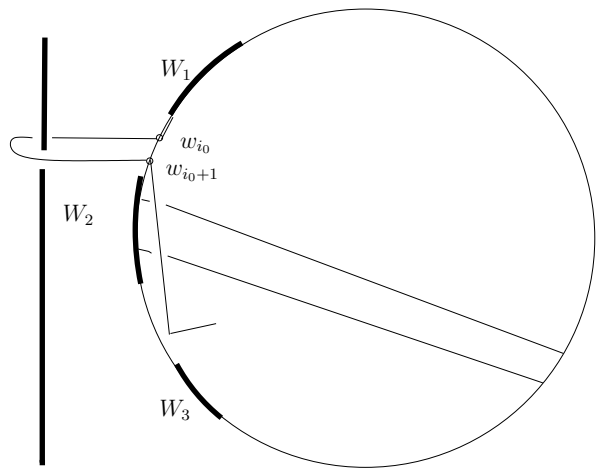}
$$
\caption{Sliding a winding, Case 2, after.}\label{figure case1b after}
\end{figure} 
\end{proof}

Proposition \ref{proposition:  layered diagrams span}  can be generalized to several other orientation schemes.  These generalizations will be useful here and in a subsequent paper ~\cite{G-cellular}.   Endow affine tangle diagrams with an arbitrary stratification.  Let $\mathcal O$  be any orientation scheme with the following properties:
\begin{enumerate}
\item  In an affine tangle diagram oriented according to $\mathcal O$, strands are oriented from lower numbered to higher numbered vertex.
\item  The order of strands is determined by some rule depending on the initial and final vertices of strands.
\item If $s$ and $t$ are two non--closed strands with initial vertices $i(s), i(t)$  and final vertices $f(s),  f(t)$, and if both $i(s) < i(t)$  and $f(s)  < f(t)$,  then $s$ precedes  $t$ in $\mathcal O$.
\item  If an affine tangle diagram $T$ is  oriented according to $\mathcal O$ and stratified,  and  if $T$ has no closed loops and no crossings of ordinary strands,  then $T$ is flagpole descending.
\end{enumerate}

Examples of such orientation schemes are the following:
\begin{itemize}
\item  Standard orientation.
\item  Final vertex  orientation:  Non--closed strands are oriented from lower numbered to higher numbered vertex and ordered according to their final vertices.
\item  Hybrid orientation:  Non--closed strands are oriented from lower numbered to higher numbered vertex and ordered as follows:    Strands with initial vertex at the top of the diagram precede those with initial vertex at the bottom of the diagram.  Strands with initial vertex at the top of the diagram are ordered according to their initial vertices.  Strands with initial vertex at the bottom of the diagram are ordered according to their final vertices.
\end{itemize}

\begin{proposition} \label{proposition:  generalization of descending diagrams span}
Endow affine tangle diagrams with an arbitrary stratification and with an orientation $\mathcal O$  satisfying the properties listed above.  Then $\akt{n, S}$ is spanned by  affine tangle diagrams  without closed loops that are
stratified according to the given stratification order and
 flagpole descending with respect to the orientation $\mathcal O$.
\end{proposition}

\begin{proof}   Let $\mathcal S$ denote the set of affine tangle diagrams without closed loops that are stratified with respect to the given stratificaton order and flagpole descending with respect to the orientation $\mathcal O$.    Let $\mathcal S'$ denote the set of affine tangle diagrams without closed loops that are stratified with respect to the given stratificaton order and flagpole descending with respect to the standard  orientation.
 
Since $\akt{n, S}$ is spanned by $\mathcal S'$ according to Proposition \ref{proposition:  layered diagrams span},  it suffices to show that $\mathcal S'$  is in the span of $\mathcal S$.  If an element of $\mathcal S'$  has no crossings of ordinary strands,  then by the assumed properties of $\mathcal O$, it is already flagpole descending with respect to $\mathcal O$.  

 Therefore, we can proceed by induction on the number of crossings.  Suppose $T \in \mathcal S'$  has $\ell \ge 1$ crossings, and assume that all elements of $\mathcal S'$  with fewer than $\ell$ crossings are in the span of $\mathcal S$.
We can assume that each strand of $T$ follows a straight line path from its initial vertex to its first crossing with the flagpole, and a straight line path from its final crossing with the flagpole to its final vertex.  (We can  change $T$ by ambient isotopy to obtain such a tangle diagram.)

If $T$ fails to be flagpole descending with respect to $\mathcal O$, then it has two strands $s$ and $t$ such that $s$ precedes $t$ in the standard orientation but  $s$ follows $t$ in the orientation $\mathcal O$, and, moreover,   such that the windings of $s$ with the flagpole lie just above the windings of $t$ with the flagpole.  (More formally:  there is a neighborhood $N$ of an interval on the flagpole such that $N$ does not intersect any strand of $T$ other than $s$ and $t$,  all windings on $s$ and $t$ lie in $N$, and the windings on $s$ lie above those on $t$.)

Notice that $i(s) < i(t)$,  since $s$ precedes $t$ in the standard orientation.   Also $f(t) < f(s)$, since otherwise,  by property (3) of $\mathcal O$,  $s$  would also precede $t$ in the orientation $\mathcal O$.
Thus we have 
$
i(s) < i(t) < f(t) < f(s).
$
The situation is illustrated in Figure \ref{figure layered 1}.

\begin{figure}[ht]
\centerline{$\inlinegraphic[scale= .9]{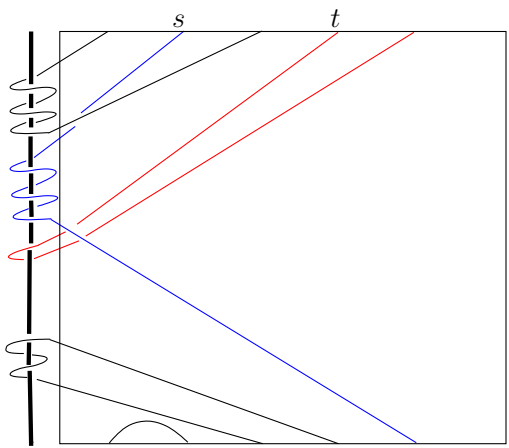}$}
\caption{Flagpole descending affine tangle diagram} \label{figure layered 1}
\end{figure}

One can change crossings on $T$ at will, because the result of changing crossings is congruent to $T$ modulo the span of affine tangle diagrams with fewer than $\ell$ crossings;  affine tangle diagrams with fewer than $\ell$ crossings are in the span the set of elements of $\mathcal S'$ with fewer than $\ell$ crossings, and these in turn are in the span of $\mathcal S$ by the induction hypothesis.

Note that $s$ and $t$ have exactly two crossings.  Changing these crossings if necessary, one can slide the last winding on $s$  below all the windings on $t$.    One can check easily that this procedure does not change the number of crossings of $s$ with $t$ or with any other strand.  Repeating this procedure several times if necessary, one can slide all windings on $s$ below the windings on $t$,  as illustrated in Figure 
\ref{figure flagpole descending 3}.  Finally,  change the crossings of $s$ and $t$ if necessary to agree with the stratification order.   The result, say $T'$,  is congruent to $T$ modulo the span of diagrams with fewer crossings,  hence modulo the span of $\mathcal S$.

\begin{figure}[h]
\centerline{$\inlinegraphic[scale= .9]{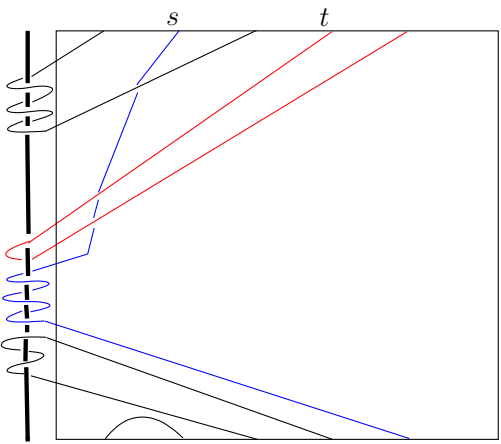}$}
\caption{Tangle diagram--after several slides} \label{figure flagpole descending 3}
\end{figure}

If $T'$ is not flagpole descending with respect to the orientation $\mathcal O$,  the procedure described in the previous paragraphs can be repeated until a flagpole descending affine tangle diagram is obtained. 
\end{proof}

We recall that a Brauer diagram is a tangle diagram in the plane,  in which information about over-- and under--crossings is ignored:

\begin{definition}\rm\label{definition brauer diagram}
An $(n,n)$--{\em Brauer diagram}  (or  $n$--{\em connector})
consists
of a collection of 
$n$ curves in the rectangle $R = I \times I$ such that
\begin{enumerate}
\item  The curves connect the points $\{\p 1,  \dots, \p n, \pbar 1,
\dots \pbar n\}$ in pairs.
\item  For each curve $C$ in the collection, the intersection of $C$ with
$\bdry(R)$ consists of the two endpoints of $C$.
\end{enumerate}
  Let $G$ be a group.   A  {\em $G$--Brauer diagram}  (or $G$--connector) is
a Brauer diagram in which each curve  (strand) is endowed with an orientation and labeled by an element of the group $G$.   Two labelings are regarded as the same if the orientation of a strand is reversed and the group element associated to the strand is inverted.
\end{definition}

We will be  interested only in $\Z$--Brauer diagrams.  
Define a map $c$ (the connector map) from oriented affine $(n,n)$--tangle diagrams  {\em without  closed loops}  to 
$\Z$--Brauer diagrams as follows.  Let $a$ be an oriented affine $(n,n)$--tangle diagram without closed loops.
 If $s$ connects two vertices
$\bm v_1$  to $\bm v_2$,   include a curve $c(s)$ in $c(a)$ connecting the same vertices with the same orientation, and label
 the oriented strand $c(s)$  with the {\em winding number} of $s$ with respect to the 
 flagpole.\footnote{The 
winding number $n(s)$ is determined 
 combinatorially as follows:  traversing the strand in its orientation, list the over--crossings $(+)$ and
under-crossings $(-)$ of the strand with the flagpole.  Cancel any
two successive  $+$'s or $-$'s in the list, so  the list now consists of
alternating $+$'s and $-$'s.  Then $n(s)$ is $\pm (1/2)$ the length of the
list, $+$ if the list begins with a $+$, and $-$ if the list begins with
a~$-$.}

\begin{lemma} \label{lemma:  isotopic layered diagrams}
Endow affine $(n,n)$--tangles diagrams with an orientation and a stratification order, each determined by some rule depending  on the vertices of the strands.
Two  affine tangle diagrams without closed loops,  with the same $\Z$--Brauer diagram, both stratified according to the given stratification order and flagpole descending with respect to the given orientation,  are ambient isotopic.
\end{lemma}

\begin{proof}    A flagpole descending affine tangle diagram without closed loops is ambient isotopic to a tight circular form,  and the tight circular form is uniquely determined by the
$\Z$--Brauer diagram of the affine tangle  diagram, the orientation, and the stratification order.
\end{proof}

\subsection{A new basis of the affine BMW algebra}
\begin{figure}[ht] 
$$
\begin{array}{c} \inlinegraphic[scale=1]{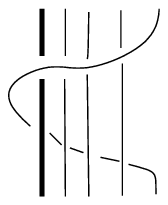} \\[-6pt]  {\scriptstyle X_4} 
\end{array}
\qquad
\begin{array}{c} \inlinegraphic{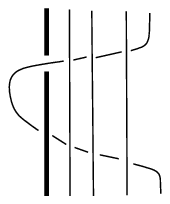} \\[-6pt]  {\scriptstyle  X_4'} 
\end{array}
$$
\caption{The elements $X_i$ and $X_i'$}
\label{figure:  X and X'}
\end{figure}

For $1 \le j \le n$,  define the affine tangle diagrams  $X_j$  and  $X_j'$ by
$$
X_j =  G_{j-1}  \cdots     G_1  X_1 G_1 \cdots G_{j-1}
$$
and
$$
X_j' =  G_{j-1}  \cdots     G_1  X_1 G_1\inv \cdots G_{j-1}\inv.
$$
See Figure \ref{figure:  X and X'}.  Also define  $Y_j' = \rho X_j' = 
G_{j-1}  \cdots     G_1  Y_1 G_1\inv \cdots G_{j-1}\inv$.

In this section we will consider affine and ordinary tangle diagrams without closed loops
 endowed with the hybrid  orientation (see the definition just before Proposition \ref{proposition:  generalization of descending diagrams span})  and the hybrid stratification order, which we now define:

 \begin{definition} \label{definition: hybrid stratification order}  
  A hybrid stratification order on the strands of an affine tangle diagram is one in which
  \begin{enumerate}
  \item  Non--closed strands precede closed loops, and strands with an initial vertex at the top of the diagram precede those with both vertices at the bottom of the diagram.
  \item Non--closed strands with initial vertex  at the top of the diagram are ordered according to the order of the initial vertices.   
  \item Non--closed strands with initial vertex at the bottom of the diagram are ordered according to the reverse of the order of the final vertices.   
 \end{enumerate}
 \end{definition}

Note that an affine tangle diagram without closed strands has a unique hybrid stratification order.

Let $d$ be a $\Z$--Brauer diagram,  and let $d_0$ be the ordinary Brauer diagram obtained by
forgetting the integer valued  labels of the strands.   
There is a unique (up to 
regular isotopy)  stratified ordinary $(n,n)$--tangle  diagram $T_{d_0}$ with no closed loops or self--crossings of strands that  has Brauer diagram $d_0$. 
Define
$$
T'_d =   ( X_{1}')^{a_1} \cdots  ( X_{n-1}')^{a_{n-1}} \  T_{d_0}  \  ( X_n')^{b_n} \cdots  ( X_1')^{b_1}, 
$$
where the exponents determined as follows:   If $d$  has a strand beginning
at a top vertex  $\p i$  with label $\ell$,  then $b_i = \ell$;  otherwise $b_i = 0$.   If $d$  has a strand beginning at a bottom vertex   and
ending at $\pbar i$ with label $\ell$,  then   $a_i = \ell$;  otherwise,  
$a_{i} = 0$.

\begin{figure}[ht]
$$
d =  \inlinegraphic{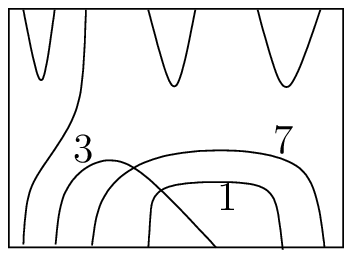}
\qquad T'_d =  \inlinegraphic{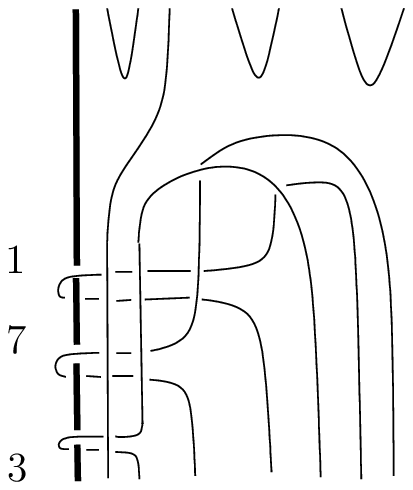}
$$
\caption{A $\Z$--Brauer diagram $d$  and its lifting $T_d'$.}
\label{figure:   Z Brauer diagram and lifting}
\end{figure}

\begin{example}  Figure \ref{figure:   Z Brauer diagram and lifting}  shows a $Z$--Brauer diagram $d$ and its lifting $T_d'$;  in the picture for $T_d'$,   the winding numbers of strands are indicated by the integers written at the left.  We have $T'_d =  ({X'}_2)^3   ({X'}_3)^7  {X'}_4  T_{d_0}$.
\end{example}

\begin{lemma} \label{lemma:  Td layered}
$T'_d$  is stratified with respect to the hybrid stratification order and flagpole descending  with respect to the hybrid orientation,  and has $\Z$--Brauer diagram equal to $d$.
\end{lemma}

\begin{proof} Straightforward.
\end{proof}

\begin{theorem} \label{theorem:  basis of affine KT}
 $\B' = \{ T'_d : d \text{ is a }  \Z  \text {--Brauer diagram}\}$  is a basis
of $\akt{n, S}$.
\end{theorem}

\begin{proof}   Let  $\rhobold$,  $\qbold$,  $\deltabold_0$,  $\deltabold_1,  \deltabold_2, \dots$
be indeterminants over $\Z$.   Let 
$$\varLambdaHat = \Z[\rhobold^{\pm 1}, \qbold^{\pm 1}, \deltabold_0^{\pm 1}, \deltabold_1, \deltabold_2, \dots]/\langle \rhobold\inv 
- \rhobold=
(\qbold\inv \qbold)(\deltabold_0 - 1) \rangle.
$$
We showed in ~\cite{GH1}, Corollary 6.14,   that for any commutative unital ring $S$ 
containing   elements 
$\rho$, $q$, and  $\delta_j$, $j \ge 0$,   with $\rho$, $q$,  and $\delta_0$ invertible, satisfying the relation
$
\rho\inv - \rho=   (q\inv -q) (\delta_0 - 1),
$
we have
$\akt{n, S}  \cong  \akt{n, \varLambdaHat} \otimes_{\varLambdaHat}   S.$ (Cf. Remark \ref{remark: change of base ring for cyclotomic BMW algebras} below.)
Therefore,  it suffices to prove the result for $\akt{n, \varLambdaHat}$.

It follows from  Proposition \ref{proposition:  generalization of descending diagrams span},  Lemma \ref{lemma:  isotopic layered diagrams}, and    Lemma \ref{lemma:  Td layered} that $\B'$  spans $\akt{n, \varLambdaHat}$.  On the other hand,
we showed in ~\cite{GH1}, Proposition 4.7,   that any collection  of affine tangle diagrams without closed loops having
distinct $\Z$--Brauer diagrams is linearly independent over $\varLambdaHat$.
\end{proof}

In $\abmw{n, S}$,  define  for $1 \le j \le n$
$$
y_j =  g_{j-1}  \cdots     g_1 \, y_1 \, g_1 \cdots g_{j-1}
$$
and
$$
y_j' =  g_{j-1}  \cdots     g_1  \, y_1  \, g_1\inv \cdots g_{j-1}\inv.
$$
For  each ordinary $n$--Brauer diagram $d_0$,  let $T_{d_0}$ be the 
unique (up to 
regular isotopy) stratified ordinary $(n,n)$--tangle  diagram with no closed loops or self--crossings of strands, and with  Brauer diagram $d_0$.   We let $\A'$ denote 
the collection of elements 
$$
   ( y_{1}')^{a_1} \cdots  ( y_{n-1}')^{a_{n-1}}  \ \varphi\inv(T_{d_0}) \   ( y_n')^{b_n} \cdots  ( y_1')^{b_1}, 
$$
where $b_i$ is zero unless $d_0$ has a strand beginning at the top vertex $\p i$   and $a_i$ is zero unless  $d_0$  has a strand  beginning at a bottom vertex   and ending at the bottom vertex $\pbar i$, and  $\varphi : \abmw{n, S} \to \akt{n, S}$  is the isomorphism of Theorem \ref{theorem:  isomorphism affine BMW and KT}.  
$\A'$ is essentially $\varphi\inv(\B')$,  with each element normalized by some power of $\rho$.

\begin{corollary}   $\A'$ is a basis of $\abmw{n, S}$.
\end{corollary}

\section{The cyclotomic BMW and Kauffman tangle algebras}
\label{The cyclotomic BMW and Kauffman tangle algebras}

\subsection{Definition of the cyclotomic BMW algebras.}
A cyclotomic  
Birman-Wenzl-Murakami algebra 
 is a quotient of the affine BMW algebra in which the affine generator $y_1$  satisfies  a monic polynomial equation.

 \begin{definition}
Let $S$ be a commutative unital ring with
parameters $\rho$, $q$, $\delta_j$   ($j \ge 0$),  and
$u_1, \dots, u_r$, with    $\rho$, $q$,  $\delta_0$, and $u_1, \dots, u_r$  invertible, and with $\rho\inv - \rho=   (q\inv -q) (\delta_0 - 1)$. 
The {\em cyclotomic BMW algebra}  $\bmw{n, S, r}(u_1, \dots, u_r)$
is the quotient of $\abmw{n, S}$ by the relation
\begin{equation} \label{equation: cyclotomic relation1}
(y_1 - u_1)(y_1 - u_2) \cdots (y_1 - u_r) = 0.
\end{equation}
\end{definition}

\begin{remark}
The assignment $e_i \mapsto e_i$, $g_i \mapsto g_i$, 
$y_1 \mapsto y_1$ defines a homomorphism $\iota$ from $\bmw  {n, S,r}$ to $\bmw  
{n+1,
S,r }$, since the relations are preserved.   It is not evident that $\iota$ is 
injective.  However,  when $S$ is an admissible integral domain (see Section \ref{section: admissiblity}),  we can show that $\bmw  {n, S,r}$ is
 isomorphic to  the cyclotomic  Kauffman tangle algebra  $\kt{n, S, r}$ (defined  below), and in this case, it is true that $\iota$ is injective.
\end{remark}

\begin{remark} \label{remark: change of base ring for cyclotomic BMW algebras}
Let $S$ be a ring with parameters $\rho$, $q$, etc., as above, and let
 $S'$  be another  ring with parameters $\rho'$, $q'$, etc.  
 Suppose there is a  ring homomorphism $\psi : S \rightarrow S'$ mapping $\rho \mapsto \rho'$,  $q \mapsto q'$,  etc.  Then $\bmw{n, S', r}$  can be regarded as an algebra over $S$, with $s x = \psi(s) x$  for $s \in S$ and $x \in \bmw{n, S', r}$, and there is an
 $S$--algebra homomorphism $\tilde\psi : \bmw{n, S, r} \to \bmw{n, S', r}$  taking generators to generators.  We claim that $\bmw{n, S, r} \otimes_S S' \cong \bmw{n, S', r}$ as $S'$--algebras.  In fact, we have the $S'$ algebra homomorphism
$\tilde\psi \otimes \id: \bmw{n, S, r} \otimes_S S'  \to \bmw{n, S', r} \otimes_S S' \cong
\bmw{n, S', r}$.  In the other direction, we have an $S'$--algebra homomorphism
$\theta: \bmw{n, S', r} \to  \bmw{n, S, r} \otimes_S S'$  mapping $g_i \mapsto g_i \otimes 1$,  etc.  The maps $\tilde \psi \otimes \id$ and $\theta$  are inverses.

Of course,  this remark applies in general to algebras defined by generators and relations.\footnote{ We could have slightly  simplified some arguments in ~\cite{GH1} using this remark. }
\end{remark}

 \subsection{Definition of the cyclotomic Kauffman tangle algebras.} 
 \label{subsection: definition of cyclotomic KT}
 Now we consider how to define a 
cyclotomic version of the Kauffman tangle algebra.
Rewrite the relation  (\ref{equation: cyclotomic relation1})  in the form
\begin{equation*} \label{equation: cyclotomic relations2} 
\sum_{k = 0}^r  (-1)^{r-k}  \varepsilon_{r-k}(u_1, \dots, u_r) y_1^k = 0,
\end{equation*}
where $\varepsilon_j$  is the $j$--th elementary symmetric function.  The corresponding relation in the affine Kauffman tangle algebra is
\begin{equation} \label{equation: cyclotomic relations3} 
\sum_{k = 0}^r  (-1)^{r-k}  \varepsilon_{r-k}(u_1, \dots, u_r) \rho^k X_1^k = 0,
\end{equation}
Now we want to impose this as a local skein relation.

 \begin{definition}
Let $S$ be a commutative unital ring with
parameters $\rho$, $q$, $\delta_j$   ($j \ge 0$),  and
$u_1, \dots, u_r$, with    $\rho$, $q$,  $\delta_0$, and $u_1, \dots, u_r$  invertible, and with $\rho\inv - \rho=   (q\inv -q) (\delta_0 - 1)$. 
The {\em cyclotomic Kauffman tangle algebra} $\kt{n, S, r}(u_1,  \dots, u_r)$ is the quotient of the affine Kauffman tangle algebra $\akt{n, S}$ by the
cyclotomic skein relation:
\begin{equation} \label{equation: kt cyclotomic relation}
\sum_{k = 0}^r  (-1)^{r-k}  \varepsilon_{r-k}(u_1, \dots, u_r) \rho^k \inlinegraphic[scale=.7]{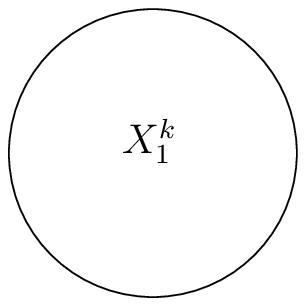} = 0,
\end{equation}
The sum is over affine tangle diagrams that are identical outside of the disc;    the interior of the disc contains an interval on the flagpole and a piece of an affine tangle diagram isotopic to $X_1^k$.
\end{definition}

We continue to write $E_i$,  $G_i$, $X_1$   for the image of these elements
of the affine Kauffman tangle algebra in the cyclotomic Kauffman tangle algebra.
We write
$$Y_1 = \rho X_1.$$

The ideal in the affine Kauffman tangle algebra by which we are taking the quotient contains the ideal generated by
$$
(Y_1 - u_1)(Y_1 - u_2) \cdots (Y_1 - u_r),
$$
but could in principal be larger.  
It follows that there is a homomorphism  $\varphi: \bmw{n, S, r}(u_1, \dots, u_r) \to \kt{n, S, r}(u_1, \dots, u_r)$  determined by $y_1 \mapsto   Y_1$,   $e_i \mapsto E_i$,  $g_i  \mapsto G_i$.  Moreover, the following diagram (in which the vertical arrows are the quotient maps) commutes.

\begin{diagram} \label{diagram: hm from cyclotomic bmw to cyclotomic kt}
\abmw {n, S }    ¤\Ear {\scriptstyle\varphi}     ¤\akt {n, S} ¤¤
\Sar {\scriptstyle\pi}         ¤          ¤\Sar {\scriptstyle\pi}  ¤¤
\movevertexleft{\bmw {n, S, r}(u_1, \dots, u_r)}                 ¤\Ear {\scriptstyle\varphi}     ¤\movevertexright{\kt {n, S, r}(u_1, \dots, u_r)}  ¤¤
\end{diagram}

\noindent
The homomorphism  $\varphi: \bmw{n, S, r}(u_1, \dots, u_r) \to \kt{n, S, r}(u_1, \dots, u_r)$  is surjective  because  the diagram commutes and
$\varphi: \abmw{n, S} \to \akt{n, S}$ is an isomorphism.

\subsection{Inclusions and conditional expectations for cyclotomic Kauffman tangle algebras.}
Let $S$  be an  ring with parameters as above and write
$\kt{n, S, r}$  for  $\kt{n, S, r}(u_1, \dots, u_r)$.
The affine Kauffman tangle algebras have inclusion maps
 $\iota : \akt {n-1, S} \rightarrow \akt {n, S}$, defined on the level of affine tangle diagrams
by adding an additional strand on the right
without adding any crossings:
$$
\iota: \quad \inlinegraphic{tangle_box} \quad \mapsto \quad 
\inlinegraphic{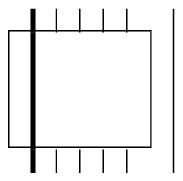}.
$$
Since these maps respect the cyclotomic relation (\ref{equation: kt cyclotomic relation}),  they induce homomorphisms
$$\iota : \kt {n-1, S, r} \rightarrow \kt {n, S, r}.$$
Recall also  that the affine Kauffman tangle algebras have a conditional expectation 
$\eps_n : \akt n \rightarrow \akt {n-1}$ defined by
$$
\eps_n(T) = \delta_0\inv{\rm cl}_n(T), 
$$
where   ${\rm cl}_n$ is the map of affine $(n,n)$--tangle diagrams to affine
$(n-1, n-1)$--tangle diagrams that ``closes" the rightmost strand, 
without adding any crossings:
$$
{\rm cl}_n: \quad \inlinegraphic{tangle_box} \quad \mapsto \quad 
\inlinegraphic{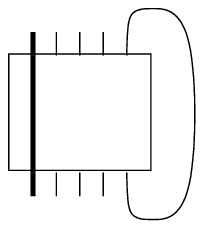}.
$$
These maps respect the cyclotomic relation (\ref{equation: kt cyclotomic relation}),  so 
induce conditional expectations
$$\eps_n : \kt{n, S, r} \rightarrow \kt{n-1, S, r}$$  
Since $\eps_n\circ \iota$ is the identity on $\kt{n-1, S, r}$,  it follows that 
$\iota : \kt {n-1, S, r} \rightarrow \kt {n, S, r}$
is injective.

\begin{remark} \label{remark: change of base ring for tangle algebras}
 Let $S$ be a ring with parameters $\rho$, $q$, etc., as above, and let
 $S'$  be another  ring with parameters $\rho'$, $q'$, etc.  
 Suppose there is a  ring homomorphism $\psi : S \rightarrow S'$ mapping $\rho \mapsto \rho'$,  $q \mapsto q'$,  etc.    Any $S'$--algebra can be regarded as an $S$--algebra using $\psi$.
 We have a map of monoid rings
$
\tilde \psi : S \ \uhat n n  \rightarrow S' \ \uhat n n ,
$
and this map respects regular isotopy,  the Kauffman skein relations, and the cyclotomic relations, so induces an $S$--algebra homomorphism from $\tilde \psi:  \kt{n, S, r} \to \kt{n, S', r}$.   As in Remark \ref{remark: change of base ring for cyclotomic BMW algebras}, we have $\kt{n,S,r} \otimes_S S' \cong \kt{n,S', r}$ as $S'$--algebras.
\end{remark}

\subsection{Finite spanning sets.}

Recall the basis $\A'$ of the affine BMW algebra $\abmw{n, S}$ described at the end of Section 2.   Let $\A'_r$  be the set of
$$
 ( y_{1}')^{a_1} \cdots  ( y_{n-1}')^{a_{n-1}} \  \varphi\inv(T_{d_0})   ( y_n')^{b_n} \cdots  ( y_1')^{b_1} \in \A' 
$$
such that $0 \le a_i, b_i  \le r-1$  for all $i$.  The cardinality of $\A'_r$ is
$r^n (2n-1)!!$.

\begin{proposition} \label{proposition: finite spanning set}
  For any  ring $S$ with with parameters $\rho$, $q$,  $\delta_j$    ($j \ge 0$), and
$u_1, \dots, u_r$, as above,  the cyclotomic BMW algebra
$\bmw{n, S, r}(u_1,  \dots, u_r)$  is spanned over $S$ by $\A'_r$.
\end{proposition}

\begin{proof}  Since $\A'$ is a basis of the affine BMW algebra
$\abmw{n, S}$,   the image of $\A'$ in $\bmw{n, S, r}(u_1, \dots, u_r)$ is spanning.
But in  the cyclotomic BMW algebra, 
each $y_j'$  is conjugate to $y_1$,   so satisfies a polynomial equation of degree $r$.
Therefore the span of $\A'_r$  equals the span of $\A'$.
\end{proof}

Let $\B'_r$  be the image of $\A'_r$ in the cyclotomic Kauffman tangle algebra, namely  $\B'_r$ is the set of 
$$
  ( Y_{1}')^{a_1} \cdots  ( Y_{n-1}')^{a_{n-1}} \  T_{d_0} \   ( Y_n')^{b_n} \cdots  ( Y_1')^{b_1} \in \B' 
$$
such that $0 \le a_i, b_i  \le r-1$  for all $i$.

 We have:

\begin{corollary} \label{corollary: finite spanning set2}
  For any  ring $S$ with with parameters $\rho$, $q$,  $\delta_j$    ($j \ge 0$), and
$u_1, \dots, u_r$, as above,  the cyclotomic Kauffman tangle algebra
$\kt{n, S, r}(u_1,  \dots, u_r)$  is spanned over $S$ by $\B'_r$.
\end{corollary}

\begin{proof} $ \varphi: \bmw {n, S, r} \to \kt {n, S, r}$  is surjective.
\end{proof}

\section{Weak admissibility and the Markov trace}
\label{section: weak admissibility and the Markov trace}

\subsection{Definition of weak admissibility} \label{subsection: weak admissibility}
The cyclotomic BMW algebras and Kauffman tangle algebras can be defined over an arbitrary commutative unital ring $S$ with parameters $\rho$, $q$,  $\delta_j$    ($j \ge 0$), and
$u_1, \dots, u_r$, such that $\rho$, $q$,  $\delta_0$  and  $u_1, \dots, u_r$   are invertible, 
and $\rho\inv - \rho=   (q\inv -q) (\delta_0 - 1)$.    However, unless the parameters
satisfy additional relations,  the identity element $\bm 1$ of the cyclotomic Kauffman tangle algebras
will be a torsion element over $S$;  if $S$ is a field (and the additional relations do not hold)  then $\bm 1 = 0$,  so $\kt {n, S, r} = \{0\}$.

Define $\delta_{-j}  \in S$  for   $j \ge 1$   by
$$
\rho^{-j} \varTheta_{-j}  = \delta_{-j} \ \bm 1.
$$
(This is an equation in the affine Kauffman tangle algebra $\akt {0, S}$).
The elements $\delta_{-j}$  can be expressed as polynomials in $\delta_0, \dots, \delta_j$  with coefficients in \break $\Z[\rho^{\pm 1}, q - q\inv]$.  In fact,  from Lemma \ref{lemma - recursion for f_r}, 
we have the recursive relations:
\begin{equation} \label{equation:  theta(-j)  recursive relations}
\begin{aligned}
\delta_{-1} &= \rho^{-2}\delta_1 \\
\delta_{-j}  & =   \rho^{ -2} \delta_j - (q\inv -q) \rho\inv \sum_{k = 1}^{j-1} (\delta_k \delta_{k-j} - \delta_{2k -j}   ).
\end{aligned}
\end{equation}
Now pass to the cyclotomic algebra $\kt {0, S, r}$.
By the cyclotomic skein relation (\ref{equation: kt cyclotomic relation}), we have for any integer $a$, 
\begin{equation*} \label{equation: kt relations among Theta curves}
\sum_{k = 0}^r  (-1)^{r-k}  \varepsilon_{r-k}(u_1, \dots, u_r) \rho^k  \varTheta_{k+a} = 0,
\end{equation*}
This gives  the relations
\begin{equation*} \label{equation: periodic relations} 
(\sum_{k = 0}^r  (-1)^{r-k}  \varepsilon_{r-k}(u_1, \dots, u_r)  \delta_{k + a})\  \bm 1  = 0,
\end{equation*}
for each $a \in \Z$.   It follows that either the identity $\bm 1$  of the cyclotomic Kauffman tangle algebra $\kt {0, S, r}$  is a torsion element over $S$, or for all $a \in \Z$, we have
\begin{equation} \label{equation: periodic relations2} 
\sum_{k = 0}^r  (-1)^{r-k}  \varepsilon_{r-k}(u_1, \dots, u_r)  \delta_{k + a}  = 0,
\end{equation}
A similar computation done in the cyclotomic BMW algebra $\bmw {2, S, r}$  shows that
$e_1$  is a torsion element unless the relations (\ref{equation: periodic relations2}) hold.

\begin{definition} 
Let $S$ be a commutative unital ring containing   elements 
$\rho$, $q$,  $\delta_j$, $j \ge 0$, and $u_1, \dots, u_r$,  with $\rho$, $q$,  $\delta_0$, and $u_1, \dots, u_r$  invertible.  We say that the parameters are {\em weakly admissible}  (or that the ring $S$ is weakly
admissible)  if the  following relations hold:
$$
\rho\inv - \rho=   (q\inv -q) (\delta_0 - 1).
$$
and 
$$
\sum_{k = 0}^r  (-1)^{r-k}  \varepsilon_{r-k}(u_1, \dots, u_r)  \delta_{k + a}  = 0,
$$
for $a \in \Z$,  where for $j \ge 1$,   $\delta_{-j}$ is defined by  the recursive relations of Equation
(\ref{equation:  theta(-j)  recursive relations}).
\end{definition}

\subsection{ $\kt {0, S, r}$ is a free $S$--module}
In this section we will show that weak admissibility of $S$ implies that $\kt {0, S, r}(u_1, \dots, u_r)$ is a free $S$--module.  A consequence of this is the existence of a special trace (the Markov trace)  on the cyclotomic Kauffman tangle algebras $\kt{n, S, r}$ and the cyclotomic BMW algebras $\bmw {n, S, r}$,  when $S$ is weakly admissible.

To show that $\kt {0, S, r}$ is a free $S$--module,  we have to show that the quotient map $\pi: \akt {0, S} \to \kt {0, S, r}(u_1, \dots, u_r)$ is injective.   We start with some observations about the kernel of $\pi: \akt {n, S} \to \kt {n, S, r}(u_1, \dots, u_r)$ for any $n$.
The kernel of $\pi$  is the span of elements $T$ of the form
$$T = \sum_{k=1}^r  (-1)^{r-k}  \varepsilon_{r-k}(u_1, \dots, u_r) \rho^k  T_k ,$$ 
where   $T_0, T_1,  \dots, T_r$  are affine $(n,n)$--tangle diagrams that are identical in the exterior of some disc $E$, and  the interior of the disc $E$ in $T_k$  contains an interval on the flagpole and a piece of an affine tangle diagram isotopic to $X_1^k$.
We refer to  such a sum as {\em generator}  of the ideal $\ker \pi$.
We can speak of the number of crossings of a generator $T$,  which is the number of crossings of ordinary strands  of  each $T_k$.

We say that $T$ is in standard position of each of its summands $T_k$ is in standard position.
The reduction by regular isotopy of an affine tangle diagram to one in standard position can be done simultaneously for the summands $T_k$ of a generator $T$,  without disturbing the interior of the disc $E$.  Therefore we can always assume that 
a generator is in standard position.

Note if we apply any of the skein relations of Definition \ref{definition affine Kauffman tangle algebra}  to each  summand $T_k$ of a generator $T$,  then we stay in the ideal $\ker \pi$.  For example,  we have
$$
T =  T' +  (q\inv - q) (T^{(0)}  - T^{(\infty)}),
$$
where $T'$ is the generator obtained by changing a certain crossing  (exterior to $E$) in each $T_k$, while     $T^{(0)}$ is the generator  obtained from the horizontal smoothing of  the same crossing in each $T_k$,  and $T^{(\infty)}$ is the generator  obtained from the vertical smoothing of  the  crossing in each $T_k$.
Thus,  $T \equiv T'$  modulo the span of  generators of $\ker \pi$  with strictly fewer crossings.

\begin{lemma}  \label{lemma:  affine 0 0 tangle with no crossings is monomial}
An affine $(0, 0)$--tangle diagram in standard position, with no crossings of ordinary strands, is equal in $\akt {0, S}$ to a monomial in $\rho\powerpm$ and 
$\{\delta_j : j \in \Z\}$.  

\end{lemma}

\begin{proof}  Let $T$ be such an affine tangle diagram.  Each closed strand comprising $T$, taken by itself,  must be flagpole--descending since it has no self--crossings.
However,  the entire affine tangle diagram is not necessarily  flagpole--descending because  simple windings on one closed strand can intervene between two simple windings on another closed strand as in the following figure:
$$
\inlinegraphic{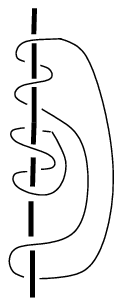}
$$

Let $A = \{p_1, \dots, p_{2k}\}$ be the set of points at which the strands of  $T$ cross the flagpole,  read from top to bottom.  Let $A_s$  be the set of points of $A$  on strand $s$.
Then $\{A_s :  s \text{ is a closed strand}\}$  is a {\em non--crossing partition} of $A$, 
see ~\cite{McCammond}.  That is,  it is not possible to have
$p_i, p_j \in A_s$  and $p_k, p_l \in A_t$  ($s \ne t$)  with $i < k < j < l$.
It follows that  some $A_s$ is an interval in $A$;   that is,  there are no
$p_i, p_j \in A_s$ and $p_k \in A_t$  ($s \ne t$)  with $i < k < j$.  But then, the strand $s$ is regularly isotopic with some $\Theta_j$  ($j \in \Z$)  and  by the free loop relations in Definition \ref{definition affine Kauffman tangle algebra}, $T = \rho^{-j}  \delta_j  \,T'$,
where $T'$ is obtained by erasing the strand $s$ from $T$.  The result now follows by induction on the number of strands of $T$.
\end{proof}

\begin{proposition}  \label{proposition:  0 cyclotomic KT algebra is free module}
Let $S$ be an  weakly admissible ring with parameters $\rho$, $q$,  $\delta_j$    ($j \ge 0$) and
$u_1, \dots, u_r$.   Then $\kt{0, S, r}(u_1,  \dots, u_r)$ is a free $S$ module of rank 1. 
\end{proposition}

\begin{proof}  We have to show that the quotient map $\pi : \akt {0, S} \to \kt {0, S, r}(u_1, \dots, u_r)$ is injective.   Let
$T = \sum_{k=1}^r  (-1)^{r-k}  \varepsilon_{r-k}(u_1, \dots, u_r) \rho^k  T_k $ 
be a generator of the ideal $\ker \pi$,  where we assume that each
$T_k$ is in standard position.

If $T$ has no crossings of ordinary strands, then by the previous lemma,
there is a monomial $m$ in $\rho\powerpm$ and $\{\delta_j : j \in \Z\}$,  and an
$a \in \Z$,   such that for each $k$,   $\rho^k T_k =  \delta_{k+a} m$.
Thus $T =  m \sum_{k = 1}^r   (-1)^{r-k}  \varepsilon_{r-k}(u_1, \dots, u_r)  \delta_{k+a}$.
But this is zero, by the weak admissibility of $S$.

We now suppose that $T$ has $\ell \ge 1$  crossings of ordinary strands,  and that
any generator $T'$ of the ideal $\ker \pi$  with fewer than $\ell$ crossings is equal to zero.  By this induction assumption  and by the remarks preceding Lemma \ref{lemma:  affine 0 0 tangle with no crossings is monomial},   changing a crossing in $T$
(that is,  changing the same crossing in each $T_k$)  does not change $T$  (as an element of $\akt {0, S}$).  So we can assume that  each $T_k$  is stratified,  in circular form with respect to a circle $D$, and tight  --- see the discussion following Remark \ref{remark:  on flagpole descending diagrams and simple windings}.

We will apply a version of the ``straightening algorithm" of Proposition \ref{proposition:  layered diagrams span}  to each $T_k$ simultaneously,  without changing anything in the interior of the disc $E$.    

Let $s_0$  be the unique strand in $T_k$ with non--trivial intersection with the disc $E$.  Suppose  $s$  is a  strand in $T_k$  other than $s_0$.   If $s$ has no  windings with the flagpole,  then $s$ can be translated away from the other strands, and then removed,  using the free loop relation.  (We have $T = \delta_0 T'$,  where $T'$ is the generator of $\ker \pi$ in which the strand $s$ has been removed from each $T_k$.)

So we suppose that $s$ does have  windings with the flagpole, and we let
 $p_1$ denote the highest winding point on  $D  \cap s$  and $p_2$  the second winding point on the same simple winding as $p_1$. Orient $s$ so that $p_1$ is the initial point and the simple winding containing $p_1$ and $p_2$ is traversed from $p_1$ to $p_2$.
Let $p_1, p_2, \dots, p_{2t}$  be the list of winding points, on $s$ and on other strands, that are below $p_1$  on $D$,  listed according to their position on $D$,  in counterclockwise order.
Let $w_1 = p_1, w_2= p_2, \dots, w_{2m}$  be the set of winding points on $s$, listed in the order  of the orientation of $s$.  We have in injection of $f: \{1, 2, \dots, 2m\} \to
\{1, 2, \dots, 2t\}$  such that $w_j = p_{f(j)}$.  Define the {\em length} of $f$ to be the number of $j$ such that $f(j) > j$;  in fact, the length depends only on the choice of $s$,  
and is the same for all $T_k$,  so we denote it by $\ell(s, T)$.

If $\ell(s, T) = 0$, then
$f(j) = j$  for all $j$,  and $s$, taken by itself,  is flagpole descending.  Moreover,  its winding points occupy an interval among all winding points on $D$, so $T_k$  is ambient isotopic to an affine tangle diagram in which $s$ is replaced by a copy of  $\Theta_{\pm m}$ having no crossings with other strands.
Then the strand $s$ can be removed,  using the free loop relation.  

Suppose now that   $\ell(s, T) > 0$.  Let $j_0$  be the first index such that $f(j_0) > j_0$.    We suppose that the simple winding with winding points $w_{j_0}$  and $w_{j_0 + 1}$  is of  type (c)  from Figure \ref{figure-simple winding types};  the other types can be handled similarly.

Now we will proceed as in the proof of Case 2 in Proposition \ref {proposition:  layered diagrams span}.
 Let $W_1$ be the set of winding points $w_i (=p_i)$ with $i < j_0$.
Let $W_2$ be the set of winding points  $p_j$   located between $W_1$  and
$w_{j_0}$  on $D$.
Let $W_3$  be the set of winding points $p_j$  located below (counterclockwise from)
   $w_{j_0 + 1}$.
Then $w_{j_0}$  is connected by a $D$--arc $\alpha$ to  the last winding point in $W_1$.  All   $D$--arcs that cross $\alpha$ are incident with $W_2$.  
 Figure  \ref{figure case1b before} illustrates the situation.  Now  exactly as in the proof of Case 2 in Proposition \ref {proposition:  layered diagrams span}, we can slide the simple winding containing the winding points $w_{j_0}$  and
 $w_{j_0 + 1}$  to a position between $W_1$ and $W_2$.  In our situation, we don't have to be concerned with any crossing changes required to do this.   We note that this move does not increase the number of crossings of $T_k$,  and does not disturb the interior of the disc $E$;  the move does decrease the length $\ell(s, T)$.   By repeating this move, we can eventually reduce the length to zero,  and then the strand $s$ can be removed by the free loop relation, as explained above.
 
 We repeat this procedure with other strands until only the strand $s_0$ having non--trivial intersection with the disc $E$ is left.  We now have
 $T = m T'$, where $m$ is a monomial in $\rho\powerpm$  and $\{\delta_j : j \in \Z\}$, 
 and $T'$ is the generator of $\ker \pi$  obtained by erasing all strands other than
 $s_0$  from all the $T_k$.     
 
 We need only a slight variation of the procedure used so far in order to deal with $s_0$.
 Let $W_0$  be the set of winding points on $D$ corresponding to simple windings in the disc $E$.  
Let $p_1, p_2, \dots, p_{2u}$  be the winding points below $W_0$ on $D$,  listed according to their position on $D$,  in counterclockwise order.   Let $q_1, q_2, \dots, q_{2t}$   be the winding points above  $W_0$ on $D$,  listed according to their position on $D$,  in counterclockwise order.   Orient $s_0$ so that the simple windings in the disc $E$  descend the flagpole.  Let $w_1, \dots, w_{2 u + 2 t}$  be the winding points on $s_0$ outside of $E$,  listed according to the orientation of $s_0$.  Define the  length of 
$T'$, denoted $\ell(T')$,    to be  $2t$ plus the number of $j$  such that $w_j = p_i$  for some $i$,  but
$j \ne i$.  If $\ell(T') = 0$, then $T'_k$ is flagpole descending and ambient isotopic to some $\Theta_{k + a}$;  it follows that $T'$ is a multiple of 
$\sum_{k = 1}^r   (-1)^{r-k}  \varepsilon_{r-k}(u_1, \dots, u_r)  \delta_{k+a}$,
which is zero, by the weak admissibility of $S$.

We suppose that $\ell(T') > 0$  and that whenever $T''$ is a generator of $\ker \pi$ such that each $T''_k$  has only one closed strand,  and $\ell(T'') <  \ell(T')$,  then
$T''= 0$.  Let $j_0$  be the first index such that $w_{j_0}  = p_i$ for some $i$ with $i \ne j_0$  or  $w_{j_0}  = q_i$ for some $i$.   If  $w_{j_0}  = p_i$ for some $i$ with $i \ne j_0$, then,  as in previously considered cases,  the length of $T'$  can be reduced by sliding a simple winding;  it follows from the induction hypothesis that $T' = 0$.  

\begin{figure}[h] \
\centerline{$\inlinegraphic[scale= .375]{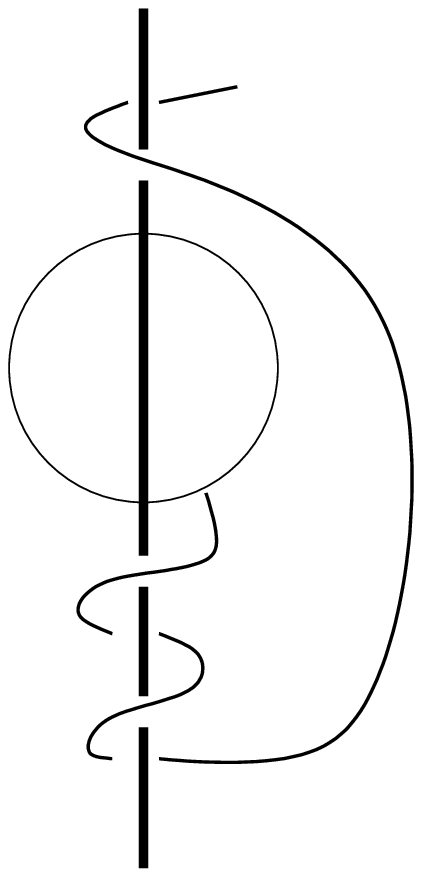}  \longrightarrow 
\inlinegraphic[scale= .375]{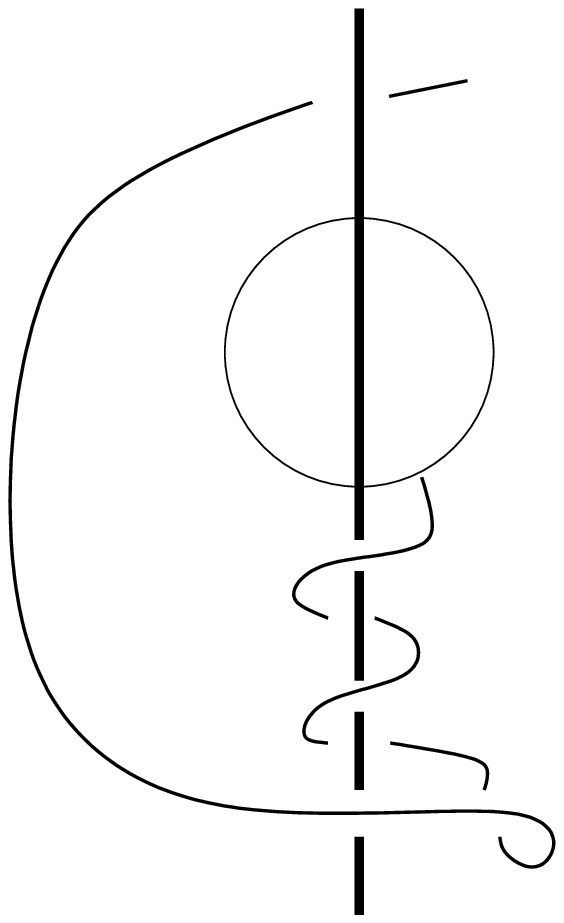}
 \longrightarrow 
\inlinegraphic[scale= .375]{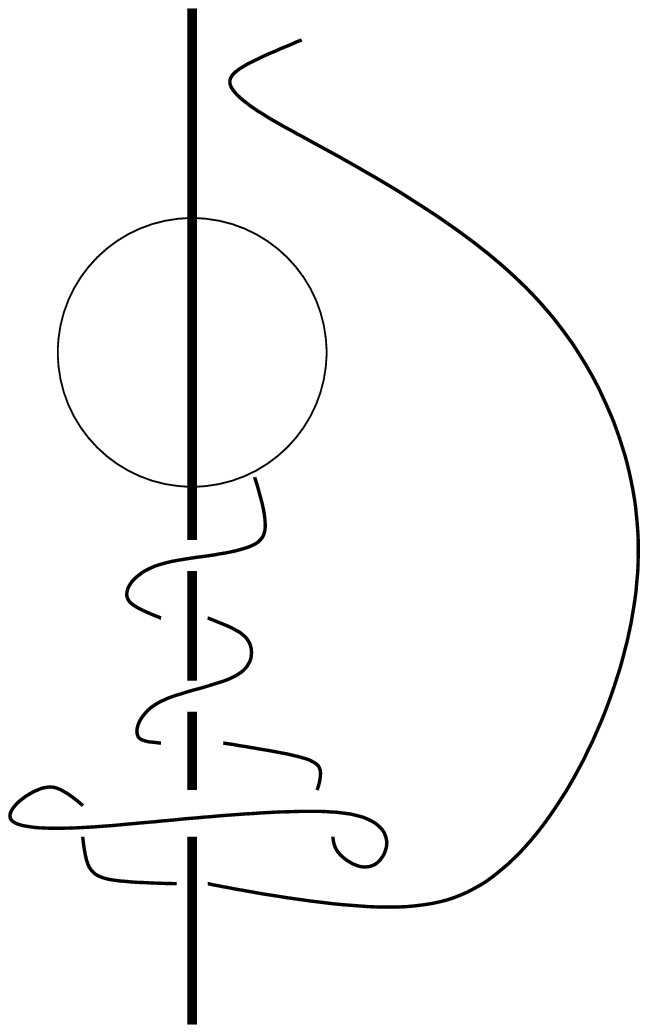}
 \longrightarrow  \rho^2
\inlinegraphic[scale= .375]{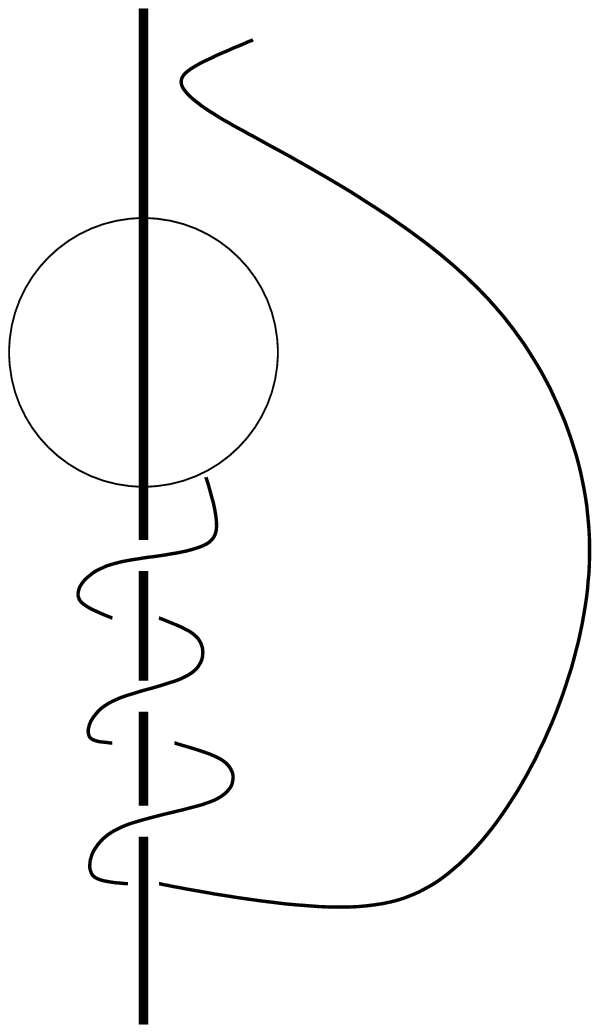}
$ 
}
\caption{A wrapping move} \label{figure a wrapping move}
\end{figure}

Suppose that
$w_{j_0} = q_i$ for some $ i$.  By pulling the arc from $w_{j_0 -1}$ to $w_{j_0}$ around the flagpole,  we can transfer the simple winding associated with the winding point 
$w_{j_0}$ from its position above the disc $E$ to a new position below  $w_{j_0 -1}$.
Suppose, for example, that the simple winding associated with the winding point 
$w_{j_0}$ is of type (d) from Figure \ref{figure-simple winding types};  then the procedure is illustrated in Figure \ref{figure a wrapping move}. Again, we do not have to be concerned with any crossing changes required to perform this move.  The move reduces the length of $T'$, so it follows from the induction hypothesis that $T' = 0$.
\end{proof}
 
\subsection{The Markov trace} 
\label{subsection: the Markov trace}
We can now identify $\kt{0, S,r}$ with $S$,  when $S$ is weakly admissible.
Recall we have conditional expectations $\eps_n: \kt{n, S, r} \to \kt{n-1, S, r}$ for 
$n \ge 1$.  It follows from ~\cite{GH1}, Proposition 2.14, that the composition
$$
\eps = \eps_{1}\circ\cdots \circ \eps_n  : \kt{n, S, r} \to \kt{0, S, r} \cong S
$$
is a trace.  We also define $\eps: \bmw{n, S, r} \to S$ by $\eps = \eps \circ \varphi$, 
where $\varphi: \bmw{n, S, r} \to \kt{n, S, r}$ is the canonical homomorphism.  Then
$\eps$ is a trace on $\bmw{n, S, R}$ with the {\em Markov property}:
 for $b \in \bmw {n-1, S, r}$,
\begin{enumerate}
\item[\rm(a)] $\eps(b g_{n-1}^{\pm 1} ) = (\rho^{\pm 1}/\delta_0) \eps(x) $,
\item[\rm(b)] $\eps(b e_{n-1}) = (1/\delta_0) \eps(x)$,
\item[\rm(c)\hskip1.2pt]  $\eps(b (y'_n)^r) = \delta_{r} \eps(b)$, and
\end{enumerate}
for $r \in \Z$, 
where $y'_n = (g_{n-1} \cdots g_1) y_1 (g_1\inv \cdots g_{n-1}\inv)$.  See ~\cite{GH1}, Corollary 6.16.

\begin{lemma} \label{lemma:  E and G not zero}
Let $S$ be a weakly admissible ring.  Then for all $n \ge 1$, $E_n$ and  $G_n$ are non-zero elements in $\kt {n, S, r}$.
\end{lemma}

\begin{proof}
$\eps(E_n) = \delta_0\inv$  and $\eps(G_n) = \rho\inv \delta_0^{-2}$.
\end{proof}

\section{Admissibility}
\label{section: admissiblity}

To obtain any substantial results about the cyclotomic BMW algebras, it appears to be necessary to impose  stronger  conditions on the ground ring.  An appropriate condition has been found by Wilcox and Yu ~\cite{Wilcox-Yu}.

Consider a commutative, unital  ring $S$ with parameters $\rho$, $q$, $\delta_j$ ($j \ge 0$)  and $u_1, \dots, u_r$, , with $\rho$, $q$, $\delta_0$,  and $u_i$ invertible,  satisfying $\rho\inv - \rho - (q\inv- q) (\delta_0 - 1)$. 
Let $a_j$ denote the signed elementary symmetric function in $u_1, \dots, u_r$, 
$$
a_j = (-1)^{r-j} \varepsilon_{r-j}(u_1, \dots, u_r).  
$$
Let
 $W_2$  denote the cyclotomic BMW algebra
$W_2 = \bmw {2, S, r}(u_1, \dots, u_r)$.  Write $e$ for $e_1$ and $g$ for $g_1$.

\begin{lemma} \label{lemma: powers of y generate W2e}    The left ideal $W_2 e$ in $W_2$  is equal to the $S$--span of  $\{e, y_1 e, \dots, y_1^{r-1} e\}$.
\end{lemma}

\begin{theorem}[Wilcox-Yu, \cite{Wilcox-Yu}] \label{theorem: equivalent conditions for admissibility}
 Let $S$ be a commutative, unital ring with
elements $\rho$, $q$, $\delta_j$ ($j \ge 0$)  and $u_1, \dots, u_r$, with $\rho$, $q$, $\delta_0$,  and $u_i$ invertible, satisfying $\rho\inv - \rho = (q\inv- q) (\delta_0 - 1)$.  
Assume that $(q - q\inv)$ is non--zero and not a zero--divisor in $S$.
The following conditions are equivalent:
\begin{enumerate}
\item  $S$ is weakly admissible, and $\{e, y_1 e, \dots, y_1^{r-1} e\}$ is linearly independent over $S$ (in $W_2 = \bmw {2, S, r}(u_1, \dots, u_r)$).
\item  The parameters satisfy the following relations:
\begin{equation} \label{equation: yu wilcox admissibility condition 1}
\begin{aligned}
&\rho(a_\ell - a_{r-\ell}/a_0) \ + 
\\& (q-q\inv)\bigg [ \sum_{j = 1}^{r - \ell} a_{j+\ell} \delta_j 
-  \sum_{j = \max(\ell + 1, \lceil r/2 \rceil)}^{\lfloor (\ell + r)/2 \rfloor} a_{2j - \ell}
+   \sum_{j =  \lceil \ell/2 \rceil}^{\min(\ell, \lceil r/2 \rceil -1)} a_{2j - \ell} \bigg ]= 0, \\& \quad \text{for $1 \le \ell \le r-1$},
\end{aligned}
\end{equation}
and
\begin{equation} \label{equation: yu wilcox admissibility condition 2}
\rho\inv a_0 - \rho a_0\inv = 
\begin{cases}
0 & \text{if  $r$ is odd} \\
(q - q\inv) & \text{if $r$ is even}.
\end{cases}
\end{equation}

\item  $S$ is weakly admissible, and $W_2 = \bmw {2, S, r}(u_1, \dots, u_r)$ admits a module $M$ with an $S$--basis $\{v_0, y_1 v_0, \dots,  y_1^{r-1} v_0\}$  such that
$e v_0 = \delta_0 v_0$.    
\end{enumerate}
\end{theorem}

\begin{definition}[Wilcox and Yu, \cite{Wilcox-Yu}]     Let $S$ be a commutative, unital ring with
elements $\rho$, $q$, $\delta_j$ ($j \ge 0$)  and $u_1, \dots, u_r$, with $\rho$, $q$, $\delta_0$,  and $u_i$ invertible, satisfying $\rho\inv - \rho = (q\inv- q) (\delta_0 - 1)$.   One says that $S$ is {\em admissible} (or that the parameters are admissible)   if $(q - q\inv)$ is non--zero and not a zero divisor in $S$ and if the equivalent conditions of Theorem \ref{theorem: equivalent conditions for admissibility} hold.
\end{definition}

It is shown in ~\cite{GH3} and in ~\cite{Yu-thesis, Wilcox-Yu2} that there exists a universal admissible integral domain $\overline R$, with the property that every admissible integral domain is a quotient of $\overline R$. Denote the parameters of 
$\overline R$ by $\rhobold$, $\qbold$,
$\deltabold_j$ ($j \ge 0$)   and $\ubold_1, \dots, \ubold_r$.   We show in ~\cite{GH3} that
$\qbold$,  $\ubold_1, \dots, \ubold_r$  are algebraically independent over $\Z$
and that  the field of fractions $F$ of  $\overline R$  is isomorphic to 
$\Q(\qbold, \ubold_1, \dots, \ubold_r)$.

The proof of the following theorem from ~\cite{GH3} depends on the existence of the Markov trace on $\kt{n, F, r}$, shown in Section \ref{section: weak admissibility and the Markov trace} of this paper, and is otherwise independent of this paper.

\begin{theorem}[\cite {GH3}]  \label{theorem:  generic structure}
Let $F$ denote the field of fractions of the universal admissible integral domain  $\overline R$.
Write $\bmw{n, F, r}$  for  $\bmw{n, F, r}(\ubold_1,  \dots, \ubold_r)$,   and
$\kt {n, F, r}$ for \break $\kt{n, F, r}(\ubold_1,  \dots, \ubold_r)$.  For all $n \ge 0$, the following assertions hold:
\begin{enumerate}
\item  $\varphi:\bmw{n, F, r} \rightarrow \kt {n, F, r}$ is an isomorphism.
\item The Markov trace $\eps$ on $\kt {n, F, r}$ is non--degenerate.
\item    $\bmw{n, F, r}$ is split semisimple of dimension  $ r^n (2n-1)!!.$
\end{enumerate}
\end{theorem}

The following theorem has been obtained independently by Wilcox and  Yu.

\begin{theorem}[Goodman and Hauschild--Mosley,  Wilcox and Yu \cite{Yu-thesis, Wilcox-Yu2}] 
\label{theorem: basis theorem over admissible integral domain}
 Let $S$ be an admissible integral domain.  Then
$\bmw{n, S, r} \cong  \kt{n, S, r}$, and $\bmw{n, S, r}$ is a free $S$--module of rank
$r^n (2n-1)!!$.
\end{theorem}

\begin{proof}  For any $S$,  the set $\A'_r$ is a spanning set in 
$\bmw{n, S, r}$ of cardinality $r^n (2n-1)!!$ by Proposition \ref{proposition: finite spanning set}.   It suffices to show that $\varphi(\A'_r) = \B'_r$ is linearly independent in 
$\kt{n, S, r}$.  When $S = F$,  this follows from Theorem \ref{theorem:  generic structure}, since $\B'_r$ is a spanning set whose cardinality equals the dimension of
$\kt{n, F, r}$.   The map $x \mapsto x \otimes 1$  from
$\kt{n, \overline R, r}$ to $\kt{n, \overline R, r} \otimes_{\overline R}  F \cong 
\kt{n, F, r}$ is $\overline R$--linear and maps  $\B'_r$ to a linearly independent set in $\kt{n, F, r}$;  hence $\B'_r$ is linearly independent in
$\kt{n, \overline R, r}$.    Finally,  since any admissible integral domain $S$ is a quotient of  $\overline R$,  and $\kt{n, \overline R, r}$ is a free $\overline R$--module with basis
$\B'_r$,   it follows that
$\kt{n, S, r} \cong \kt{n, \overline R, r} \otimes_{\overline R} S$  is a free $S$--module
with basis $\B'_r$.
\end{proof}.

\section{Remarks on the affine and cyclotomic Hecke algebras}

In this section, we apply our techniques to the affine and cyclotomic Hecke algebras, rather than the affine and cyclotomic BMW algebras,  recovering the main technical result (Theorem 5) of Lambropoulou ~\cite{LambrJKTR1999}.  This result was used in ~\cite{LambrJKTR1999} to construct Markov traces on the Artin braid group of type $B$ that factor through the affine Hecke algebra, and thus invariants of links in the solid torus.

First we want to explain that the Artin braid group of type $B$ is isomorphic to the group of braids in the annulus cross the interval.  Consequently,  the affine Hecke algebra  can be identified with the algebra of such braids, modulo Hecke skein relations.  This is proved,
for example, in  ~\cite{crisp}  and ~\cite{allcock-braids},  but we want to point out  an elementary proof, using only facts from Section 1.4 of ~\cite{Birman-book},  and a short argument from the proof of 
~\cite{crisp-paris},  Proposition 2.1.  (We suppose that this elementary proof must also be well known.)

A {\em braid}  is an ordinary tangle,  in the disk cross the interval, in which each strand is monotone; that is, each strand intersects every horizontal plane exactly once.  We identify ambient isotopic braids.  Braids can be represented by braid diagrams, 
which are  ordinary tangle diagrams in which each strand is monotone.
  The set of braids with a given number of  strands forms a group under composition of tangles.    The $n$--strand braid group will be denoted by 
$\mathcal B_n$.   We remind the reader that our convention for the composition $a b$
of braid diagrams is that $b$ is stacked over $a$.

Let $\sigma_i$ denote the braid diagram$$\sigma_i = \inlinegraphic{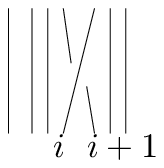}.$$   Artin showed that $\braid n$ has a presentation with generators $\sigma_1, \dots, \sigma_{n-1}$  and relations
$\sigma_i \sigma_{i+1}  \sigma_i = \sigma_{i+1}  \sigma_i \sigma_i $,  and $\sigma_i \sigma_j = \sigma_j \sigma_i$ if $|i - j| \ge 2$.

An {\em affine braid}  is a tangle in the annulus cross the interval, in which each strand is monotone;  again, we identify ambient isotopic affine braids.   Affine braids can be represented by affine braid diagrams,  which are affine tangle diagrams in which each strand is monotone.
The set of affine braids with a given number strands forms a group under composition of affine tangles.
We denote the group of affine braids with $n$  strands by $\affbraid n$.

Every braid $b$ with $n$ strands induces a permutation $\pi(b) \in S_n$.   If $b$ has a strand connecting the $k$--th upper vertex with the $j$--th lower vertex, then  
$\pi(b)(k) = j$.    The pure braid group $\mathcal P_n$  is the subgroup of braids inducing the identity permutation.    

Fix an integer $n$ for the remainder of this section.  

 Let $\braid {n+1}$  be the braid group on $n+1$ strands with vertices labelled by $0, 1, \dots, n$.    The affine braid group $\affbraid n$ can be identified with the set of those braids in  $\braid {n+1}$ having a strand connecting the $0$--th upper vertex with the 
$0$--th lower vertex. Let  $\braid n$   denote the subgroup generated by $\sigma_1, \dots, \sigma_{n-1}$.

\begin{definition}  The Artin group  $\mathcal A (B_n)$   of type $B_n$  is the group with generators 
\break $\beta_0, \beta_1, \dots, \beta_{n-1}$,   and defining relations
\begin{enumerate}
\item  $\beta_0 \beta_1 \beta_0 \beta_1 =  \beta_1 \beta_0 \beta_1 \beta_0$.
\item $\beta_i \beta_{i+1}  \beta_i = \beta_{i+1}  \beta_i \beta_i $  for $1 \le i \le n-2$.
\item $\beta_i \beta_j = \beta_j \beta_i$ if $|i - j| \ge 2$.
\end{enumerate}
\end{definition}

It is easy to check that $\beta_0 \mapsto \sigma_0^2$,  and $\beta_i \mapsto \sigma_i$ for $1 \le i \le n-1$  defines a homomorphism  $\varphi: \mathcal A(B_n) \rightarrow \affbraid n  \subset \braid {n+1}$.   We are going to show that this map is an isomorphism.

\begin{proposition}[\cite{crisp}, \cite{allcock-braids}] \label{proposition: artin group and affine braid group}
 The Artin group $\mathcal A(B_n)$  is isomorphic to the affine braid group $\affbraid n$.
\end{proposition}

\def\s{\sigma}
\begin{proof}
For $0 \le i < j \le n$,  let $$A_{i, j} = \s_{j-1} \cdots \s_{i+1}\  \s_i^2 \   \s_{i+1}\inv \cdots \s_{j-1}\inv.$$
By computation, or by a picture proof, one verifies that
$$
A_{i, j} =  \s_i\inv \cdots \s_{j-2}\inv \ \s_{j-1}^2 \ \s_{j-2} \cdots \s_i.
$$
The following  is part of the statement of   ~\cite{Birman-book},  Lemma 1.8.2.

\noindent {\bf Fact 1:}   The set of $A_{i, j}$  generates the pure braid group $\mathcal P_{n+1}$.

Let $V_0$  denote the subgroup of $\mathcal P_{n+1}$  generated by the set of 
$A_{0, j}$  with $1 \le j \le n$.  One can check that $V_0$  is normalized by 
$\s_j$  for $j \ge 1$.  In fact, we have

\begin{subequations}  \label{equation:  A_ij relations}
\begin{align}
\s_j A_{0,j} \s_j\inv &=  A_{0, j+1} \quad (1 \le j \le n-1) \label{equation:  1st A_ij relation}\\
\s_k A_{0,j} \s_k\inv  &= A_{0,j} \quad (k \ge 1 \text{ and } k  \not\in \{ j-1, j\}), 
\label{equation:  2nd A_ij relation}  \\
\s_j A_{0,j+1} \s_j\inv &=  A_{0,j+1}\inv A_{0,j} A_{0,j+1} \quad (1 \le j \le n-1).
\label{equation:  last A_ij relation}
\end{align}
\end{subequations}

The first two relations  are checked easily by picture proofs or by using the Artin presentation.   The last relation can proved by induction on $j$.  The base case
$j = 1$ follows from the braid relations.   Assume relation (\ref{equation:  last A_ij relation}) for a particular value of $j$  and apply $\Ad( \s_j \s_{j+1})$ to both sides of the equation;   reducing using the braid relations  as well as relations (\ref{equation:  1st A_ij relation}) and
(\ref{equation:  2nd A_ij relation})   yields  (\ref{equation:  last A_ij relation})  with $j+1$ in place of $j$.

It follows from Fact 1 that $V_0$ is normal in $\mathcal P_{n+1}$, and
$ \mathcal P_{n+1} = \mathcal P_{n} V_0$.
  Moreover,  $V_0$ is normalized by $\braid n$.   
  
  Now any element  $b \in \braid {n+1}$  can be written as
$ b = \gamma b_0$,  where $b_0$ is a pure braid and $\gamma$ is a permutation braid,
i.e. a braid in which any two strands cross at most once.  Moreover, $\pi(b) = \pi(\gamma)$.   If $b \in \affbraid n$,  then $\gamma \in \braid n$.  That is, 
$$
\affbraid n \subseteq \braid n  \mathcal P_{n+1} =  \braid n \mathcal P_n  V_0 = \braid n V_0  \subseteq \affbraid n.
$$
Thus we have

\noindent {\bf Fact 2:}   \qquad $ \affbraid n = \braid n \ltimes V_0.$

It follows that the homomorphism  $\varphi: \mathcal A(B_n) \rightarrow \affbraid n  $ is surjective.

The map $\s_i \mapsto \s_{n-1-i}\inv$  determines an automorphism of $\braid {n+1}$ which takes $$A_{0,j} = \s_0\inv \cdots \s_{j-2}\inv\  \s_{j-1}^2 \s_{j-2} \cdots \s_0$$ to
$$\s_{n-1} \cdots \s_{n-j+1} \  \s_{n-j}^{-2} \  \s_{n-j+1}\inv \cdots \s_{n-1}\inv = A_{n-j, n}\inv.$$
The image of $V_0$ is the group $U_n$  generated by $\{A_{k, n} :  0 \le k \le n-1\}$.
By ~\cite{Birman-book},  page 23,   $U_n$ is a free group with free basis $\{A_{k, n} :  0 \le k \le n-1\}$.  Hence,

\noindent {\bf Fact 3:}   $V_0$ is a free group with free basis $\{A_{0, j} : 1 \le j \le n\}$.

Now we can attempt to define a homomorphism $\psi : \affbraid n \rightarrow \mathcal A(B_n)$  via the requirements  $\psi(\s_i) = \beta_i$  for $1 \le i \le n-1$   and
$$\psi(A_{0, j}) =  t_j :=  \beta_{j-1} \cdots \beta_1 \ \beta_0 \  \beta_1\inv \cdots \beta_{j-1}\inv.$$
Because the $\beta_i$  for $i \ge 1$   satisfy the ordinary braid relations,  and since
the set of $A_{0, j}$  are free generators of $V_0$, these requirements define
homomorphisms on $\braid n$  and on $V_0$.  To check that $\psi$ extends to a homomorphism on   $ \affbraid n = \braid n \ltimes V_0$,  it suffices to check that
$\beta_i t_j \beta_i\inv =  \psi(\s_i  A_{0, j} \s_i \inv)$  for $i, j \ge 1$.  That is,  we have to check that

\begin{subequations}  \label{equation:  t_j relations}
\begin{align}
\beta_j t_j \beta_j\inv &=  t_{j+1} \quad (1 \le j \le n-1) \label{equation:  1st t_j relation}\\
\beta_k t_j \beta_k\inv  &= t_j \quad (k \ge 1 \text{ and } k  \not\in \{ j-1, j\}), 
\label{equation:  2nd t_j relation}  \\
\beta_j t_{j+1} \beta_j\inv &=  t_{j+1}\inv t_j t_{j+1} \quad (1 \le j \le n-1).
\label{equation:  last t_j relation}
\end{align}
\end{subequations}
But these relations can be checked in exactly the same way as the relations 
(\ref{equation: A_ij relations}).

Finally,  we have $\psi\circ\varphi(\beta_j)  =  \beta_j$  for all $j$,  so $\psi\circ\varphi$ is the identity on $\mathcal A(B_n)$.
\end{proof}

\begin{definition}
Let $S$ be a commutative unital ring with an invertible element $q$.  
\begin{enumerate}
\item  The ordinary Hecke algebra $\hec {n,S}(q^2)$ of type $A$ is the quotient of the group algebra  $S\, \braid n$
of the braid group,  by the relations  $$\s_i - \s_i\inv =  (q - q\inv) \quad  (1 \le i \le n-1).$$
\item  The affine Hecke algebra $\ahec {n, S}(q^2)$  is the quotient of
the group algebra   $S\, \mathcal A(B_n)$ of the Artin group of type $B_n$,  by the relations 
$$\beta_i - \beta_i\inv =  (q - q\inv) \quad  (1 \le i \le n-1).$$
\end{enumerate}
\end{definition}

Since $\braid n$ imbeds in $\affbraid n$,  the ordinary Hecke algebra imbeds in the affine Hecke algebra.  For $i \ge 1$,   we denote the image of $\s_i$ in the ordinary Hecke algebra  (and the image of $\beta_i$ in the affine Hecke algebra)  by $g_i$.
We denote the image of $\beta_0$ in the affine Hecke algebra  by $x_1$ and
define
$$
x_j =  g_{j-1}  \cdots     g_1  x_1 g_1 \cdots g_{j-1},
$$
and
$$
x_j' =  g_{j-1}  \cdots     g_1  x_1 g_1\inv \cdots g_{j-1}\inv,
$$
for $1 \le j \le n$.

\begin{corollary} Let $S$ be any commutative ring with identity and with an invertible element $q$.  The affine Hecke algebra  $\ahec {n,S}(q^2)$   is isomorphic to the $S$--algebra of affine braid diagrams,  modulo the Hecke skein relation:
$$
\qquad\quad \inlinegraphic[scale=.7]{pos_crossing} - \inlinegraphic[scale=.35]{neg_crossing} 
\quad = 
\quad
(q - q\inv)\,
\inlinegraphic[scale=2]{id_smoothing}.
$$
\end{corollary}
Here, the figures indicate affine braid diagrams which  are identical outside the region shown.

\begin{proof}   By Proposition ~\ref{proposition: artin group and affine braid group},
the Artin group  $\mathcal A(B_n)$ is isomorphic to the group of affine braid diagrams.
Hence the affine Hecke algebra is isomorphic to the algebra of affine braid diagrams, 
modulo the ideal generated by the relations $\s_i - \s_i\inv = (q - q\inv)$.  Since any affine braid diagram is isotopic to a product of the elementary diagrams $\s_i$,  $x_1$, and their inverses,  the ideal  in $S \affbraid n$  generated by the Hecke skein relations is the same as the ideal generated by the relations $\s_i - \s_i\inv = (q - q\inv)$.
\end{proof}

It is well known that the affine Hecke algebra $\ahec {n, S}(q^2)$  is a free module over the ordinary Hecke algebra $\hec {n, S}(q^2)$ , with basis consisting of Laurent monomials in the commuting elements $x_j$.  
We want to use the technique from Section 2 of this paper
 to show that the affine Hecke algebra has a basis as an
$\hec {n,S}(q^2)$  module  consisting of ordered  Laurent monomials in the non--commuting elements $x'_j$:
$$
(x'_n)^{a_n}  \cdots  (x'_1)^{a_1}
$$ 
This is a theorem of Lambropoulou (\cite{LambrJKTR1999}, Theorem 5),  and is the Hecke algebra analogue of our
Theorem \ref{theorem:  basis of affine KT}.

When an affine braid diagram is written as a word in the generators $\s_i^{\pm 1}$ and
$x_1^{\pm 1}$,  it is already in standard position as an affine tangle diagram, cf.  Definition \ref{definition: standard position}.    We give affine braid diagrams the standard orientation,  so strands are oriented downward,  and ordered according to the order of their initial vertices, from left to right.   We use the same order for the stratification order,  so an affine braid diagram is stratified if it is totally descending;  that is,  each crossing is encountered first as an over crossing.

\begin{lemma}  $\ahec {n, S}(q^2)$ is spanned by totally descending affine braid diagrams.
\end{lemma}

\begin{proof}  Same as the proof of Lemma \ref{lemma-totally descending tangles
span}.
\end{proof}

\begin{lemma} \label{lemma:  totally descending flagpole descending braids span}
 $\ahec {n, S}(q^2)$ is spanned by totally descending and flagpole descending affine braid diagrams.
\end{lemma}

\begin{proof}  Same as the proof of Proposition \ref{proposition:  layered diagrams span}
(except there is no need to deal with closed loops).
\end{proof}

The connector map $c$  (see the paragraph before Lemma \ref{lemma:  isotopic layered diagrams})  maps  affine braid diagrams to $\Z$--permutation diagrams, i.e., $\Z$--Brauer diagrams with vertical strands only.  The set of  $\Z$--permutation diagrams  constitutes a multiplicative group in the $\Z$--Brauer algebra, isomorphic to the wreath product  $\Z \wr S_n$.

\begin{lemma}  \label{lemma:  isotopic affine braid diagrams}

Two totally descending, flagpole descending affine braid diagrams with the same $\Z$--permutation diagram are isotopic.
\end{lemma}

\begin{proof}  This is a special case of Lemma \ref{lemma:  isotopic layered diagrams}.
\end{proof}

Let $d$ be a $\Z$--permutation diagram.   Let $\alpha$  be the underlying permutation, and let $a_j$  be the $\Z$--valued label of the $j$--th strand.   Let $T_\alpha$ be the unique totally descending ordinary tangle diagram with underlying permutation 
$\alpha$.   (If $s_{i_1} \cdots s_{i_r}$ is a reduced expression for $\alpha$, then
$T_\alpha =   g_{i_1}\inv \cdots  g_{i_r}\inv$.)    Define
$$
T'_d =  T_\alpha \  (x'_n)^{a_n}\cdots(x'_1)^{a_1}
$$
It is straightforward to check that $T'_d$  is totally descending,  flagpole descending and has $\Z$--permutation diagram equal to $d$.

\begin{proposition}[\cite{LambrJKTR1999}, Theorem 5]  
\label{proposition:  x' basis in affine Hecke algebra}
For any commutative unital ring $S$,  
$$\Sigma' = \{ T'_d :  d \text{ is a $\Z$--permutation diagram}\}$$   is an $S$--basis of $\ahec {n, S}(q^2)$.
\end{proposition}

\begin{proof}  Consider the generic ground ring for the affine Hecke algebra,
$A = \Z[\qbold, \qbold\inv]$,   where $\qbold$ is an indeterminant.    Since the affine Hecke algebra over $A$ is a free $A$--module,  it follows that for any ring $S$, 
$$\ahec {n, S}(q^2)  =  \ahec {n, A}(\qbold^2) \otimes_{A} S.$$  Therefore it suffices to prove the result for $\ahec {n, A}(\qbold^2) $.

By Lemmas \ref{lemma:  totally descending flagpole descending braids span}  and \ref{lemma:  isotopic affine braid diagrams},   $\Sigma'$
spans $\ahec {n, A}(\qbold^2) $. It remains to show that $\Sigma'$  is linearly independent over $A$.   Suppose we have a linear relation:
$
\sum_d  r_d  T'_d = 0
$
in  $\ahec {n, A}(\qbold^2) $.   We can suppose that the non--zero coefficients are polynomials in $\Z[\qbold]$  and that they have no common factor,  hence no common integer root.   

The specialization  $\ahec {n, \Q}(1)$  with $q = 1$  can be identified with the group algebra of $\Z \wr S_n$,  and the element $T'_d$ is thus identified with the group element $d$.   These elements are linearly independent over $\Q$,  so the relation
$\sum_d r_d(1) T'_d = 0$  implies that $r_d(1) = 0$ for all $d$.   Since the non--zero
$r_d$ have no common integer root,  we must have $r_d = 0$ for all $d$.
\end{proof}

\begin{remark}  The passage from Theorem 4 to Theorem 5 in ~\cite{LambrJKTR1999} is reversible, i.e.,  the two theorems provide two different descriptions of the same basis of $\ahec {n, S}(q^2)$.    Therefore,  Proposition \ref{proposition:  x' basis in affine Hecke algebra} implies the existence of Markov traces on the affine Hecke algebra and of Jones type invariants of links in the solid torus,  as in  ~\cite{LambrJKTR1999}, Sections 4 and 5.
\end{remark}

\begin{definition}  Let $S$ be a commutative unital ring with invertible element $q$.
Let $r \ge 1$, and let
 $u_1,  \dots, u_r$  be additional invertible elements in $S$.   The {\em cyclotomic Hecke algebra}
$\hec{n, S, r}(q^2; u_1,  \dots, u_r)$  is the quotient of the affine Hecke algebra $\ahec{n, S}(q^2)$ by the polynomial relation    $(x_1 - u_1) \cdots (x_1 - u_r) = 0$.
\end{definition}

Let $\Sigma'_r$   denote the set of all  $ T_\alpha   \  (x'_n)^{a_n}\cdots(x'_1)^{a_1} \in \Sigma'$  such that $  0 \le a_j \le r-1 $  for all $j$.

\begin{corollary}  For any commutative unital ring $S$ with invertible elements
$q, u_1,  \dots, u_r$,  the cyclotomic Hecke algebra  $\hec{n, S, r}(q^2; u_1,  \dots, u_r)$ is a free $S$--module with basis   $\Sigma'_r$.
\end{corollary}

\begin{proof}  Let $R = \Z[\qbold^{\pm 1},  \boldu_1^{\pm 1},  \dots,  \boldu_r^{\pm 1}]$.
For any $S$,
$$
\hec{n, S, r}(q^2; u_1,  \dots, u_r) =  \hec{n, R, r}(\qbold^2; \boldu_1,  \dots, \boldu_r) \otimes_R S.
$$
Therefore, it suffices to prove the result for the cyclotomic Hecke algebra over $R$.

As in the proof of Proposition \ref{proposition: finite spanning set},  it follows from Proposition \ref{proposition:  x' basis in affine Hecke algebra} that
$\Sigma'_r$  spans  $ \hec{n, R, r}(\qbold^2; \boldu_1,  \dots, \boldu_r)$.  
Moreover,  $\Sigma'_r$  is linearly independent over the field of fractions 
$F$  of $R$, because the cyclotomic Hecke algebra over $F$ is an $F$--vector space of dimension $r^n n!$,  and $\Sigma'_r$ is a spanning set of the same cardinality.
\end{proof}

\bibliographystyle{amsplain}
\bibliography{cyclotomicBMW}

\end{document}